\documentclass[10.5pt]{article}
\usepackage{bbm}
\usepackage[leqno]{amsmath}
\usepackage{amsfonts}
\usepackage{graphicx}

\usepackage{amsmath}
\usepackage{amssymb}
\usepackage{latexsym}
\usepackage{amsmath, amsfonts,amssymb, amsthm, euscript,makeidx,color,mathrsfs}

\usepackage{enumerate} 

\def\sqr#1#2{{\vcenter{\vbox{\hrule height.#2pt
              \hbox{\vrule width.#2pt height#1pt \kern#1pt \vrule width.#2pt}
              \hrule height.#2pt}}}}
\def\signed #1{{\unskip\nobreak\hfil\penalty50
              \hskip2em\hbox{}\nobreak\hfil#1
              \parfillskip=0pt \finalhyphendemerits=0 \par}}
\def\endpf{\signed {$\sqr69$}}

\def\3n{\negthinspace \negthinspace \negthinspace }
\def\2n{\negthinspace \negthinspace }
\def\1n{\negthinspace }

\def\dbE{\mathbb{E}}
\def\dbF{\mathbb{F}}

\def\dbH{\mathbb{H}}

\def\dbP{\mathbb{P}}

\def\dbR{\mathbb{R}}
\def\dbS{\mathbb{S}}

\def\sD{\mathscr{D}}

\def\sU{\mathscr{U}}

\def\={\buildrel \triangle \over =}

\def\ds{\displaystyle}

\def\ns{\noalign{\ss}}
\def\a{\alpha}

\def\d{\delta}

\def\l{\lambda}
\def\m{\mu}
\def\n{\nu}
\def\si{\sigma}
\def\t{\tau}
\def\f{\varphi}
\def\th{\theta}

\def\G{\Gamma}
\def\D{\Delta}
\def\Th{\Theta}

\def\O{\Omega}

\def\cF{{\cal F}}

\def\cK{{\cal K}}

\def\cP{{\cal P}}

\def\ss{\smallskip}
\def\ms{\medskip}

\def\q{\quad}
\def\qq{\qquad}
\def\hb{\hbox}

\def\lan{\mathop{\langle}}
\def\ran{\mathop{\rangle}}

\def\essinf{\mathop{\rm essinf}}

\def\wt{\widetilde}

\def\cd{\cdot}
\def\cds{\cdots}

\def\as{\hbox{\rm a.s.{ }}}

\def\tr{\hbox{\rm tr$\,$}}

\def\les{\leqslant}
\def\ges{\geqslant}

\def\({\Big (}
\def\){\Big )}
\def\[{\Big[}
\def\]{\Big]}
\def\bde{\begin{definition}}
\def\ede{\end{definition}}
\def\be{\begin{equation}}
\def\bel{\begin{equation}\label}
\def\ee{\end{equation}}
\def\bt{\begin{theorem}}
\def\et{\end{theorem}}
\def\bc{\begin{corollary}}
\def\ec{\end{corollary}}
\def\bl{\begin{lemma}}
\def\el{\end{lemma}}
\def\bp{\begin{proposition}}
\def\ep{\end{proposition}}
\def\bas{\begin{assumption}}
\def\eas{\end{assumption}}
\def\br{\begin{remark}}
\def\er{\end{remark}}
\def\ba{\begin{array}}
\def\ea{\end{array}}
\def\ed{\end{document}}

\def\square#1{\vbox{\hrule\hbox{\vrule height#1%
     \kern#1\vrule}\hrule}}
\def\rectangle#1#2{\vbox{\hrule\hbox{\vrule height#1%
     \kern#2\vrule}\hrule}}

\font\tenbb=msbm10 \font\sevenbb=msbm7 \font\fivebb=msbm5

\newfam\bbfam
\scriptscriptfont\bbfam=\fivebb \textfont\bbfam=\tenbb
\scriptfont\bbfam=\sevenbb

\newtheorem{lemma}{Lemma}[section]
\newtheorem{remark}{Remark}[section]

\newtheorem{theorem}{Theorem}[section]
\newtheorem{corollary}{Corollary}[section]

\newtheorem{definition}{Definition}[section]
\newtheorem{proposition}{Proposition}[section]
\newtheorem{assumption}{Assumption}[section]

\makeatletter
   
   \@addtoreset{equation}{section}
\makeatother

\begin{document}
\title{\bf Time-Inconsistent Recursive Stochastic Optimal Control Problems\thanks{This work is supported in part
by NSF Grant DMS-1406776, the National Natural Science Foundation of
China (11471192, 11401091,11571203), and the Nature Science Foundation of Shandong
Province (JQ201401).}}

\author{Qingmeng Wei\footnote{School of Mathematics and Statistics, Northeast Normal University, Changchun 130024, China; email: {\tt weiqm100@nenu.} {\tt edu.cn}},~~
Jiongmin Yong\footnote{Department of Mathematics, University of
Central Florida, Orlando, FL 32816, USA; email: {\tt
jiongmin.yong@ucf.edu}},~~and~~Zhiyong Yu\footnote{Corresponding
author, School of Mathematics, Shandong University, Jinan 250100,
China; email: {\tt yuzhiyong@sdu.edu.cn}}}

\maketitle

\begin{abstract}
In this paper, we study a time-inconsistent stochastic optimal control problem with a recursive cost functional by a multi-person hierarchical differential game approach. An equilibrium strategy of this problem is constructed and a corresponding equilibrium Hamilton-Jacobi-Bellman (HJB, for short) equation is established to characterize the associated equilibrium value function. Moreover, a well-posedness result of the equilibrium HJB equation is established under certain conditions.
\end{abstract}


\ms


\section{Introduction}\label{Sec_Introduction}

Let $(\O,\cF,\dbF,\dbP)$ be a complete filtered
probability space on which a $d$-dimensional standard Brownian
motion $W(\cd)$ is defined, whose natural filtration is $\dbF
=\{\cF_t\}_{t\ges0}$ (augmented by all the $\dbP$-null
sets). Let $T>0$. We denote
$$\sD=\Big\{(t,\xi)\bigm| t\in[0,T],\hb{ and $\xi$ is
$\cF_t$-measurable, $\dbR^n$-valued, with }\dbE
|\xi|^2<\infty\Big\}.$$
For any $(t,\xi)\in\sD$, we consider the following controlled
stochastic differential equation (SDE, for short):
\bel{Sec1_Sys}\left\{\2n\ba{ll}
\ns\ds dX(s)=b(s,X(s),u(s))ds+\si(s,X(s),u(s))dW(s),\qq s\in[t,T],\\
\ns\ds X(t)=\xi,\ea\right.\ee
where $b:[0,T]\times\dbR^n\times U\to\dbR^n$ and
$\si:[0,T]\times\dbR^n\times U\to\dbR^{n\times
d}$ are suitable deterministic maps with $U\subseteq\dbR^m$
being a nonempty set, which could be bounded or unbounded. In the
above, $(t,\xi)\in\sD$ is called an {\it initial pair},
$u:[t,T]\times\O\to U$ is called a {\it control
process}, and $X:[t,T]\times\O\to\dbR^n$ is called
a {\it state process}. We define the set of all {\it admissible
control processes} by the following
$$\sU[t,T]=\Big\{u:[t,T]\times\O\to U\bigm|u(\cd)\hb{ is $\dbF$-progressively measurable, with }
\dbE\int_t^T|u(s)|^2ds <\infty\Big\}.$$
Under some mild conditions, for any $(t,\xi)\in\sD$ and $u(\cd)\in\sU[t,T]$, \eqref{Sec1_Sys} admits a unique strong solution $X(\cd)\equiv X(\cd\,;t,\xi,u(\cd))$. To measure the performance of the control process $u(\cd)\in\sU[t,T]$, we may introduce the following cost functional
\bel{Sec1_Cost_Cons}J^0(t,\xi;u(\cd))=\dbE_t\[e^{-\l(T-t)}h^0(X(T))
+\int_t^T e^{-\l(s-t)}g^0(s,X(s),u(s))ds\],\ee
with some constant $\l\ges0$ (called a {\it discount rate}), some maps $g^0:[0,T]\times\dbR^n\times U\to\dbR$ and $h^0:\dbR^n\to\dbR$, and $\dbE_t[\,\cd\,]=\dbE[\,\cd\,|\cF_t]$. On the right hand side of
\eqref{Sec1_Cost_Cons}, the first term is referred to as a discounted {\it terminal cost}, and the second term is referred to as a discounted {\it running cost}. We note that discount terms
$e^{-\l(T-t)}$ and $e^{-\l(s-t)}$ appearing in the terminal cost and running cost are exponential functions which have the same parameter $\l$. We therefore call \eqref{Sec1_Cost_Cons} a cost functional
with an {\it exponential discount}. Now, let us state the following classical optimal control problem.

\ms

{\bf Problem (C)$^0$.} For any $(t,\xi)\in\sD$, find a $\bar u(\cd)\in\sU[t,T]$ such that
\bel{Sec1_Optimal}J^0(t,\xi;\bar u(\cd))=\essinf_{u(\cd)\in\sU[t,T]}J^0(t,\xi;u(\cd)).\ee

If an admissible control $\bar u(\cd)$ satisfies \eqref{Sec1_Optimal}, it is called an {\it optimal control} for Problem (C)$^0$ with the initial pair $(t,\xi)$. The corresponding
state process $\bar X(\cd)$ and $(\bar X(\cd),\bar u(\cd))$ are called an {\it optimal state}
and an {\it optimal pair}. A function $V:[0,T]\times\dbR^n\to\dbR$ is called a
{\it value function} of Problem (C)$^0$ if
\bel{}V(t,x)=J(t,x;\bar u(t)),\q\dbP\,\hb{-}\as, \
\forall\,(t,x)\in[0,T]\times\dbR^n.\ee

For the above problem, it is standard that (see, for example, \cite{Yong2011}) if $\bar
u(\cd)$ is an optimal control for the initial pair $(t,\xi)\in\sD$ with the corresponding
{\it optimal state process} $\bar X(\cd)$, then for any $s\in[t,T]$,
$$J^0(s,\bar X(s);\bar u(\cd)|_{[s,T]})=\essinf_{u(\cd)\in\sU[s,T]}J^0(s,\bar X(s);u(\cd)).$$
This means that the restriction $\bar u(\cd)|_{[s,T]}$ of $\bar
u(\cd)$ on $[s,T]$ is an optimal control for the corresponding
initial pair $(s,\bar X(s))$. Such a property is referred to as the
{\it time-consistency} of the optimal control $\bar u(\cdot)$, or
the time-consistency of Problem (C)$^0$.

\ms

If we let $(Y_0(\cd),Z_0(\cd))\equiv(Y_0(\cd\,;t,\xi,u(\cd)),Z_0(\cd\,;t,\xi,u(\cd)))$ be the
adapted solution to the following (linear) {\it backward stochastic differential equation} (BSDE, for short):
\bel{Sec1_BSDE_Cons}\left\{\2n\ba{ll}
\ns\ds dY_0(s)=-\[\l Y_0(s)+g^0(s,X(s),u(s))\]ds+Z_0(s)dW(s),\qq s\in[t,T],\\
\ns\ds Y_0(T)=h^0(X(T)),\ea\right.\ee
then the following holds:
$$Y_0(t)=\dbE_t\bigg[e^{-\l(T-t)}h^0(X(T))+\int_t^T
e^{-\l(s-t)}g^0(s,X(s),u(s))ds\bigg]=J^0(t,\xi;u(\cd)).$$
Therefore, $Y_0(t)\equiv Y_0(t;t,\xi,u(\cd))$ is a representation
of $J^0(t,\xi,u(\cd))$. Further, from BSDE \eqref{Sec1_BSDE_Cons},
we also have
$$Y_0(t)=\dbE_t\[h^0(X(T))+\int_t^T\(\l Y_0(s)+g^0(s,X(s),u(s))\)ds\].$$
This reminds us the {\it stochastic differential utility} (SDU, for
short) introduced by Duffie and Epstein in 1992
(\cite{Duffie-Epstein1992a,Duffie-Epstein1992b}, see also
\cite{Duffie-Lions1992}), which is the continuous-time version of {\it
recursive utility} (\cite{Kraft-Seifried2014}). More precisely, for
a terminal random payoff $\eta$ (at time $T$) and a consumption
process $c(\cd)$ on $[t,T]$ (which can be regarded as a control
process), a corresponding SDU process $Y(\cd)$ is a solution to
the following equation:
\bel{Sec1_SDU}
Y(s)=\dbE_s\[\eta+\int_s^Tg(r,c(r),Y(r))dr\],\qq s\in [t,T],\ee
for some proper map $g(\cd)$. It is by now well-understood that
under some mild conditions, the following (possibly non-linear) BSDE
admits a unique adapted solution $(Y(\cd),Z(\cd))$:
$$\left\{\2n\ba{ll}
\ns\ds dY(s)=-g(s,c(s),Y(s))ds+Z(s)dW(s),\qq s\in[t,T],\\
\ns\ds Y(T)=\eta,\ea\right.$$
and the process $Y(\cdot)$ satisfies \eqref{Sec1_SDU}. In another
word, a general SDU process can be represented by the adapted
solution to a proper BSDE. Because of this, the SDU of
Duffie-Epstein has been extended substantially later (\cite{El
Karoui-Peng-Quenez1997,Antonelli-Barucci-Mancino2001,Lazrak-Quenez2003,Lazrak2004}).
Hereafter, we call the adapted solution $(Y(\cdot),Z(\cdot))$ of a
BSDE an SDU process.

Inspired by the SDUs, for control system \eqref{Sec1_Sys}, we may
introduce the following cost functional:
\bel{Sec1_Cost_Cons_Rec}J(t,\xi;u(\cd))=Y(t;t,\xi,u(\cd)),\ee
with $(Y(\cd),Z(\cd))\equiv(Y(\cd\,;t,\xi,u(\cd)),Z(\cd\,;t,\xi,u(\cd)))$ being the adapted solution to the following BSDE:
\bel{Sec1_BSDE_Cons_Rec}\left\{\2n\ba{ll}
\ns\ds dY(s)=-g(s,X(s),u(s),Y(s),Z(s))ds+Z(s)dW(s),\qq s\in[t,T],\\
\ns\ds Y(T)=h(X(T)),\ea\right.\ee
for which \eqref{Sec1_BSDE_Cons} is a special case.

\ms

Based on the above, we rewrite \eqref{Sec1_Sys} and
\eqref{Sec1_BSDE_Cons_Rec} in a compact form:
\bel{}\left\{\2n\ba{ll}
\ns\ds dX(s)=b(s,X(s),u(s))ds+\si(s,X(s),u(s))dW(s),\qq s\in[t,T],\\
\ns\ds dY(s)=-g(s,X(s),u(s),Y(s),Z(s))ds+Z(s)dW(s),\qq s\in[t,T],\\
\ns\ds X(t)=\xi,\qq Y(T)=h(X(T)),\ea\right.\ee
which is a controlled decoupled {\it forward-backward stochastic differential
equation} (FBSDE, for short, see \cite{Ma-Yong1999} for some
relevant results), and consider the cost functional defined by
\eqref{Sec1_Cost_Cons_Rec}. With these, we may formally pose the
following optimal control problem involving SDU:

\ms

{\bf Problem (C).} For any $(t,\xi)\in\sD$, find a $\bar
u(\cd)\in\sU[t,T]$ such that
$$J(t,\xi;\bar u(\cd))=\essinf_{u(\cd)\in\sU[t,T]}J(t,\xi;u(\cd)).$$

Similar to Problem (C)$^0$, one can show that Problem (C) is also time-consistent
(\cite{Peng1997}).

\ms

Although time-consistency is a very good feature of Problems (C)$^0$
and (C), it is a little too ideal. In real world, time-consistent
situation rarely exists. Instead, most problems, if not all, people
encounter are not time-consistent. In another word, an optimal
policy/control found for the current initial pair $(t,\xi)$ will
hardly stay optimal as time goes by. We refer to such kind of
situations the {\it time-inconsistency}. Some careful observations
show that one important reason leading to time-inconsistency is due
to people's subjective {\it time-preferences}. As a matter of fact,
people usually discount more on the utility for the outcome of
immediate future events. Mathematically, such a situation can be
described by the so-called {\it non-exponential discounting},
meaning that the discounting terms $e^{-\l(T-t)}$ and
$e^{-\l(s-t)}$ appear in the terminal cost and the running cost
are replaced by some more general functions $\m(T,t)$ and
$\n(s,t)$, respectively. We note that even if these functions are
respectively replaced by exponential functions $e^{-\l_1(T-t)}$
and $e^{-\l_2(s-t)}$, as long as $\l_1\ne\l_2$, the
corresponding problem will already become time-inconsistent
(\cite{MN2010,MS2011}). As suggested in
\cite{Yong2011,Yong2012a,Yong2012b,Yong2013,Yong2014}, instead of
\eqref{Sec1_Cost_Cons}, one may consider cost functional
\bel{Sec1_Cost_Inc}J(t,\xi;u(\cd))=\dbE_t\[h(t,X(T))+\int_t^T
g(t,s,X(s),u(s))ds\].\ee
It is clear that the classical situation \eqref{Sec1_Cost_Cons}
corresponds to the following special case:
$$h(t,x)=e^{-\l(T-t)}h^0(x),\qq g(t,s,x,u)=e^{-\l(s-t)}g^0(s,x,u).$$
However, we see easily that the above \eqref{Sec1_Cost_Inc} does not
contain the problems involving SDUs. In order to include problems
involving SDUs, we now propose the following:
\bel{Sec1_BSDE_Inc_Rec}\left\{\ba{ll}
\ns\ds dY(s)=-g(t,s,X(s),u(s),Y(s),Z(s))ds+Z(s)dW(s),\qq s\in[t,T],\\
\ns\ds Y(T)=h(t,X(T)),\ea\right.\ee
then we may let the cost functional to be
\bel{Sec1_Cost_Inc_Rec}J(t,\xi;u(\cd))=Y(t;t,\xi,u(\cd))=\dbE_t\[
h(t,X(T))+\int_t^Tg(t,s,X(s),u(s),Y(s),Z(s))ds\].\ee

Let us make a couple of observations. Firstly, if
$$g(t,s,x,u,y,z)\equiv g(t,s,x,u),$$
then \eqref{Sec1_Cost_Inc_Rec} is reduced to \eqref{Sec1_Cost_Inc},
and \eqref{Sec1_BSDE_Inc_Rec} is not necessary. Secondly, if
$$h(t,x)\equiv h(x),\qq g(t,s,x,u,y,z)\equiv g(s,x,u,y,z),$$
then \eqref{Sec1_BSDE_Inc_Rec} is reduced to
\eqref{Sec1_BSDE_Cons_Rec} and \eqref{Sec1_Cost_Inc_Rec} coincides
with \eqref{Sec1_Cost_Cons_Rec}, so that the corresponding optimal
control problem becomes Problem (C).

\ms

Suggested by the above, we may now introduce the following controlled decoupled FBSDE:
\bel{Sec1_FBSDE_Inc}\left\{\2n\ba{ll}
\ns\ds dX(s)=b(s,X(s),u(s))ds+\si(s,X(s),u(s))dW(s),\qq s\in[t,T],\\
\ns\ds dY(s)=-g(t,s,X(s),u(s),Y(s),Z(s))ds+Z(s)dW(s),\qq s\in[t,T],\\
\ns\ds X(t)=\xi,\qq Y(T)=h(t,X(T)).\ea\right.\ee
Under some mild conditions, for any $(t,\xi)\in\sD$ and $u(\cd)\in\sU[t,T]$, the above
FBSDE admits a unique adapted solution $(X(\cd),Y(\cd),Z(\cd))$ (\cite{Ma-Yong1999}).
Then we may define the recursive cost functional by
\bel{recursive}J(t,\xi;u(\cd))=Y(t)\equiv Y(t;t,\xi,u(\cd)),\ee
and an optimal control problem can be posed. It is expected that such an optimal control problem is time-inconsistent. Therefore, finding an optimal control at any given initial pair $(t,\xi)$ is not very useful. Instead, one should wisely find an equilibrium strategy which is time consistent and possesses certain kind of local optimality.

\ms

To find time-consistent equilibrium strategy, we adopt the method of multi-person differential game. The idea can be at least traced back to the work of Pollak \cite{Pollak1968} in 1968. Later, the approach was adopted and further developed by Ekeland--lazrak \cite{Ekeland-Lazrak2010}, Yong \cite{Yong2012b,Yong2014,Yong2016}, Bj\"{o}rk--Murgoci \cite{Bjork-Murgoci2014}, and Bj\"{o}rk--Murgoci--Zhou \cite{Bjork-Murgoci-Zhou2014} for various kinds of problems. Let us now elaborate the approach a little more carefully to our stochastic recursive cost case as follows, which has some substantial and interesting differences from the works mentioned above.

\ms

Firstly, we divide the whole time interval $[0,T]$ into $N$ subintervals: $[t_0,t_1),[t_1,t_2),\cds,[t_{N-1},t_N]$, with $t_0=0$, $t_N=T$, and introduce an $N$-person differential game, where players are labeled from 1 through $N$. Player $k$ takes over the system at time $t_{k-1}$ from Player $(k-1)$, and controls the system on $[t_{k-1},t_k)$, then hand it over to Player $(k+1)$ at $t_k$. The ``sophisticated'' recursive cost functional for Player $k$ is defined through a BSDE on $[t_{k-1},t_k]$, whose coefficient/generator depends on his/her initial pair $(t_{k-1}, X^k(t_{k-1}))$ with $X^k(t_{k-1})$ equal to $X^{k-1}(t_{k-1})$ (the terminal state of Player $(k-1)$) and whose terminal value at $t_k$ equals $\Th^k(t_k,X^k(t_k))$ with $X^k(t_k)$ being the terminal state of Player $k$. The function $\Th^k(\cd\,,\cd)$ is constructed based on the assumption that later players will play optimally with respect to their ``sophisticated'' recursive cost functionals. Therefore, the ``sophisticated'' recursive cost functionals are constructed {\it recursively}. On the other hand, although he/she will not control the system on $[t_k,T]$, Player $k$ will still ``discounts'' the future costs in his/her own way, due to the time-preference feature of the problem \cite{Yong2012b,Yong2014,Yong2016}. It turns out that each player faces a resulted time-consistent optimal control problem. Therefore, under suitable conditions, each player will have an optimal control defined on the corresponding subinterval. Then we could construct a partition-dependent equilibrium strategy and the corresponding partition-dependent equilibrium value function of the game.

\ms

Secondly, letting the mesh size of the partition tend to zero, we (at least formally could) get the limits called the {\it time-consistent equilibrium strategy} and {\it time-consistent equilibrium value function} of the original time-inconsistent optimal control problem. At the same time, a so-called {\it equilibrium Hamilton--Jacobi--Bellman equation} (equilibrium  HJB equation, for short) is also derived to characterize the time-consistent equilibrium value function. Moreover, in the case that the equilibrium HJB equation is well-posed, the formal convergence (as the mesh size goes to zero) of the relevant functions will become rigorous.

\ms

Finally, at the moment, to establish the well-posedness of the equilibrium HJB equation, we will assume that the diffusion term $\si$ of the state equation does not depend on the control process $u(\cd)$, beside some other mild conditions. The general case that $\si$ contains the control process $u(\cd)$ is still open and will be investigated in the future.

\ms

By the way, as we mentioned, the optimal control problem for every player (with the ``sophisticated'' recursive cost functional) is time-consistent. Thus, it is expected that one could use the classical approach to deal with them. To this end, we establish a stochastic verification theorem. To our best knowledge, there is no existing ready-to-use result which can be applied to our problem directly.

\ms

The rest of this paper is organized as follows. Section \ref{Sec_Pre} is devoted to the preliminaries for our study. We recall the relationship between FBSDEs and PDEs, and establish a verification theorem for time-consistent optimal control problem with recursive cost functional. In Section \ref{Sec_Formulation}, the time-inconsistent recursive stochastic optimal control problem is formulated. Then we introduce and solve the multi-person differential game in Section \ref{Sec_Game}, which leads to an approximate time-consistent equilibrium strategy for the original problem. By letting the mesh size of the partition go to zero, we also formally obtain an equilibrium strategy and a characterization of the equilibrium value function in terms of the equilibrium HJB equation. Finally, the well-posedness of the equilibrium HJB equation is established under proper conditions.

\section{Preliminaries}\label{Sec_Pre}

We introduce the following notation:
$$D[0,T]=\Big\{(t,s)\in[0,T]^2\bigm|0\les t\les s\les T\Big\}.$$
Let $U\subseteq\dbR^m$ be a nonempty set, which could be bounded or unbounded. Let maps $b:[0,T]\times\dbR^n\times U\to\dbR^n$, $\si:[0,T]\times\dbR^n\times U\to\dbR^{n\times d}$, $g:D[0,T]\times\dbR^n\times U\times\dbR\times\dbR^{1\times d}\to\dbR$, and $h:[0,T]\times\dbR^n\to\dbR$ satisfy the following assumptions.

\ms

{\bf (H1)} Maps $b$, $\si$, $g$ and $h$ are continuous and there exists a constant $L>0$ such that, for any $(t,s,u)\in D[0,T]\times U$, any $u_1,u_2\in U$, $x_1,x_2\in\dbR^n$, $y_1,y_2\in\dbR$, $z_1,z_2\in\dbR^{1\times d}$,
$$\ba{ll}
\ns\ds|b(s,x_1,u)-b(s,x_2,u)|+|\si(s,x_1,u)-\si(s,x_2,u)|\\
\ns\ds\qq+|g(t,s,x_1,u,y_1,z_1)-g(t,s,x_2,u,y_2,z_2)|+|h(t,x_1)-h(t,x_2)|\\
\ns\ds\les L(|x_1-x_2|+|y_1-y_2|+|z_1-z_2|),\ea$$
and
$$|b(s,0,u)|+|\si(s,0,u)|+|g(t,s,0,u,0,0)|\les L(1+|u|).$$

\ms

The following standard result shows that FBSDE
\eqref{Sec1_FBSDE_Inc} is well-posed.

\begin{proposition}\label{Sec2_Prop_WellPosedness_FBSDE}
\sl Let {\rm(H1)} hold. Then for any $(t,\xi)\in\sD$ and
$u(\cd)\in\sU[t,T]$, FBSDE \eqref{Sec1_FBSDE_Inc} admits a
unique adapted solution $(X(\cd),Y(\cd),Z(\cd))\equiv
(X(\cd\,;t,\xi,u(\cd)),Y(\cd\,;t,\xi,u(\cd)),Z(\cd\,;t,\xi,u(\cd)))$
such that
$$\dbE_t\[\sup_{t\les s\les T}|X(s)|^2+\sup_{t\les s\les T}|Y(s)|^2+\int_t^T|Z(s)|^2ds\]\les K\Big\{1+|\xi|^2+\dbE_t\int_t^T|u(s)|^2 ds\Big\}.$$

\end{proposition}

Before going further, we introduce some notations. Let $\dbS^n\subseteq\dbR^{n\times n}$ be the set of all
$(n\times n)$ symmetric matrices. Let
\bel{Sec2_Notations}\left\{\2n\ba{ll}
\ns\ds a(t,x,u)={1\over2}\si(t,x,u)\si(t,x,u)^\top,\qq(t,x,u)\in[0,T]\times\dbR^n\times U,\\
\ns\ds\dbH(\t,t,x,u,\th,p,P)=\tr\big[a(t,x,u)P\big]+\big\langle b(t,x,u),\ p\big\rangle
+g\big(\t,t,x,u,\th,p^\top\si(t,x,u)\big),\\
\ns\ds\hskip 5cm (\t,t,x,u,\th,p,P)\in D[0,T]\times\dbR^n\times U\times\dbR\times\dbR^n\times\dbS^n,\ea
\right.\ee
where the superscript $\top$ denotes the transpose of vectors or matrices. We note that, since $U$ is not necessarily compact, $\ds\inf_{u\in U}\dbH(\t,t,x,u,\th,p,P)$ may be infinite on the whole space $D[0,T]\times\dbR^n\times U\times\dbR\times\dbR^n\times\dbS^n$. Similar to Yong
\cite{Yong2012b}, we introduce the following assumption.

\ms

{\bf (H2)} There exists a map $\psi:D[0,T]\times\dbR^n\times\dbR\times\dbR^n\times\dbS^n\to U$
with needed regularity such that
$$\ba{ll}
\ns\ds\psi(\t,t,x,\th,p,P)\in\,\hb{argmin }\dbH(\t,t,x,\cd\,,\th,p,P)\\
\ns\ds\qq\qq\qq\q\equiv\Big\{\bar u\in U\bigm|\dbH(\t,t,x,\bar
u,\th,p,P)=\min_{u\in U}\dbH(\t,t,x,u,\th,p,P)\Big\},\\
\ns\ds\qq\qq\qq\qq\qq\qq(\t,t,x,\th,p,P)\in D[0,T]\times\dbR^n\times\dbR\times\dbR^n\times\dbS^n.\ea$$
For more explanations and comments on the above Assumption (H2), one
is referred to \cite{Yong2012b}.

\subsection{FBSDEs and PDEs}

As preparations, we begin with the following family of FBSDEs
without involving controls, which are time-consistent and
parameterized by the initial pairs $(t,\xi)\in\sD$:
\bel{Sec2.1_FBSDE}\left\{\2n\ba{ll}
\ns\ds dX(s)=b(s,X(s))ds+\si(s,X(s))dW(s),\qq s\in[t,T],\\
\ns\ds dY(s)=-g(s,X(s),Y(s),Z(s))ds+Z(s)dW(s),\qq s\in[t,T],\\
\ns\ds X(t)=\xi,\qq Y(T)=h(X(T)).\ea\right.\ee
Note that, under Assumption (H1) (ignoring $t$ and $u$), by
Proposition \ref{Sec2_Prop_WellPosedness_FBSDE},
\eqref{Sec2.1_FBSDE} admits a unique adapted solution
$(X(\cd),Y(\cd),Z(\cd))\equiv
(X(\cd\,;t,\xi),Y(\cd\,;t,\xi),Z(\cd\,;t,\xi))$. We point out that
Assumption (H1) can be substantially relaxed still guaranteeing the
existence and uniqueness of the adapted solution to FBSDE \eqref{Sec2.1_FBSDE}. As
suggested in \cite{Peng1991,Pardoux-Peng1992,Ma-Protter-Yong1994},
the family of FBSDEs \eqref{Sec2.1_FBSDE} is closely linked to the
following semi-linear partial differential equation (PDE, for
short):
\bel{Sec2.1_PDE}\left\{\2n\ba{ll}
\ns\ds\Th_t(t,x)+\dbH\big(t,x,\Th(t,x),\Th_x(t,x),\Th_{xx}(t,x)\big)=0,\qq
(t,x)\in [0,T]\times\dbR^n,\\
\ns\ds\Th(T,x)=h(x),\qq x\in\dbR^n,\ea\right.\ee
where, for simplicity, we use the notation $\dbH$ defined in
\eqref{Sec2_Notations} omitting $\t$ and $u$.

We denote
$$C^{1,2}([0,T]\times\dbR^n)=\Big\{v(\cd\,,\cd)\in
C([0,T]\times\dbR^n)\bigm|v_t(\cd,\cd),v_x(\cd,\cd), v_{xx}(\cd,\cd)\in C([0,T]\times\dbR^n)
\Big\}.$$
The following result will be used below.

\begin{theorem}\label{Sec2.1_THM_Feynman-Kac}
\sl Suppose $\Th(\cd\,,\cd)\in C^{1,2}([0,T]\times\dbR^n)$ is
a classical solution to PDE \eqref{Sec2.1_PDE}. For any given
$(t,\xi)\in\sD$, suppose FBSDE \eqref{Sec2.1_FBSDE} admits a
unique solution $(X(\cd),Y(\cd),Z(\cd))\equiv
(X(\cd\,;t,\xi),Y(\cd\,;t,\xi),Z(\cd\,;t,\xi))$. Then
$$\Th(t,\xi)=Y(t;t,\xi),\qq\as$$

\end{theorem}

\it Proof. \rm For any given initial pair $(t,\xi)\in\sD$, noting that
$X(\cdot)$ is the solution to the forward equation in \eqref{Sec2.1_FBSDE},
by applying  It\^o's formula to $\Th(\cd\,,X(\cd))$, we get
$$\ba{ll}
\ns\ds\Th(s,X(s))=h(X(T))-\int_s^T\Big\{\Th_t(r,X(r))+\tr
\big[a(r,X(r))\Th_{xx}(r,X(r))\big]+\lan b(r,X(r)),\Th_x(r,X(r))\ran\Big\}dr\\
\ns\ds\qq\qq\qq\qq\qq-\int_t^T\Th_x(r,X(r))^\top\si(r,X(r))dW(r),\qq s\in[t,T].\ea$$
By the definition of function $\dbH$ (see \eqref{Sec2_Notations} ignoring $\t$ and $u$), and noticing
$\Th(\cd\,,\cd)$ satisfies PDE \eqref{Sec2.1_PDE}, the above equation can be rewritten as
$$\ba{ll}
\ns\ds\Th(s,X(s))=h(X(T))-\int_s^Tg\big(r,X(r),\Th(r,X(r)),\Th_x(r,X(r))^\top\si(r,X(r))\big)dr\\
\ns\ds\qq\qq\qq\qq-\int_s^T\Th_x(r,X(r))^\top\si(r,X(r))dW(r),\qq s\in[t,T].\ea$$
By the uniqueness of the backward equation in \eqref{Sec2.1_FBSDE}, we have
$$\Th(s,X(s))=Y(s;t,\xi),\qq\as,\forall s\in[t,T].$$
Consequently, by letting $s=t$, we obtain our conclusion. \endpf

\subsection{Problem (C) and verification theorem}

In this subsection, we consider the controlled form of \eqref{Sec2.1_FBSDE}, i.e., for any initial pair $(t,\xi)\in\sD$, the controlled FBSDE is given by
\bel{Sec2.2_FBSDE}\left\{\2n\ba{ll}
\ns\ds dX(s)=b(s,X(s),u(s))ds+\si(s,X(s),u(s))dW(s),\qq s\in[t,T],\\
\ns\ds dY(s)=-g(s,X(s),u(s),Y(s),Z(s))ds+Z(s)dW(s),\qq s\in[t,T],\\
\ns\ds X(t)=\xi,\qq Y(T)=h(X(T)).\ea\right.\ee
For any initial pair $(t,\xi)\in\sD$ and any control process
$u(\cd)\in\sU[t,T]$, under Assumption (H1) ignoring $t$,
Proposition \ref{Sec2_Prop_WellPosedness_FBSDE} works again to
ensure the existence and uniqueness of solution to
\eqref{Sec2.2_FBSDE}. As in Section \ref{Sec_Introduction}, now we
introduce the following recursive cost functional
\bel{}
J(t,\xi;u(\cd)=Y(t;t,\xi,u(\cd)),\qq u(\cd)\in\sU[t,T],\ee
and present a result about Problem (C).

\ms

Comparing with PDE \eqref{Sec2.1_PDE} which is related to the
situation without involving control processes, the family of
controlled FBSDEs \eqref{Sec2.2_FBSDE} are closely linked to the
following so-called Hamilton-Jacobi-Bellman (HJB, for short)
equation which is a fully non-linear PDE:
\bel{Sec2.2_HJB}\left\{\2n\ba{ll}
\ns\ds V_t(t,x)+\inf_{u\in U}\dbH\big(t,x,u,V(t,x),V_x(t,x),V_{xx}(t,x)\big)=0,\qq
(t,x)\in [0,T]\times\dbR^n,\\
\ns\ds V(T,x)=h(x),\qq x\in\dbR^n,\ea\right.\ee
where, for simplicity, we use the notation $\dbH$ defined in
\eqref{Sec2_Notations} ignoring $\t$ once again.

\ms

The following result is called a verification theorem for Problem
(C) which can be regarded as a generalization of Theorem
\ref{Sec2.1_THM_Feynman-Kac}.

\begin{theorem}\label{Sec2.2_THM_Verification}
\sl Let {\rm(H1)} hold. Suppose $V(\cd\,,\cd)\in C^{1,2}([0,T]\times\dbR^n)$ is a classical solution
to the HJB equation \eqref{Sec2.2_HJB}. Then, for any $(t,\xi)\in\sD$,
\bel{Sec2.2_VT_1}
V(t,\xi)\les J(t,\xi;u(\cd)),\qq\as,\ \forall\
u(\cd)\in\sU[t,T].\ee
Moreover, let $(\bar X(\cd),\bar u(\cd))$ be a state-control
pair with $(t,\xi)$ such that
\bel{Sec2.2_VT_2}
\bar u(s)\in\hb{\rm argmin }\dbH(s,\bar X(s),\cd\,,V(s,\bar
X(s)),V_x(s,\bar X(s)),V_{xx}(s,\bar X(s))),\qq s\in [t,T],\ee
Then
\bel{Sec2.2_VT_3}V(t,\xi)=J(t,\xi;\bar u(\cd)),\qq\as\ee
In another word, $V(\cd\,,\cd)$ is the value function of Problem
(C), and $(\bar X(\cd),\bar u(\cd))$ is an optimal pair of
Problem (C) for the initial pair $(t,\xi)$.
\end{theorem}

\it Proof. \rm For any $(t,\xi)\in\sD$, let $(X(\cd),u(\cd))$ be an
admissible pair with $(t,\xi)$. From the HJB equation \eqref{Sec2.2_HJB},
\bel{Sec2.2_VT_Temp1}V_t(s,X(s))+\dbH\big(s,X(s),u(s),V(s,X(s)),V_x(s,X(s)),V_{xx}(s,X(s))\big)\ges 0,\qq
s\in[t,T].\ee
On the other hand, by applying It\^{o}'s formula to $V(\cd\,,X(\cd))$, we have
$$\ba{ll}
\ns\ds V(s,X(s))=h(X(T))-\int_s^T\Big\{V_t(r,X(r))+\tr\big[
a(r,X(r),u(r))V_{xx}(r,X(r))\big]\\
\ns\ds\qq\qq\qq\qq\qq\qq\qq+\big\langle b(r,X(r),u(r)),\ V_x(r,X(r))\big\rangle\bigg\}dr\\
\ns\ds\qq\qq\qq\qq\qq\qq-\int_s^T V_x(r,X(r))^\top\si(r,X(r),u(r))dW(r),\qq s\in[t,T].\ea$$
By the definition of function $\dbH$ (see \eqref{Sec2_Notations} ignoring $\t$), the above equation can be
rewritten as
\bel{Sec2.2_VT_Temp2}\ba{ll}
\ns\ds V(s,X(s))=h(X(T))+\int_s^T\Big\{g\big(r,X(r),u(r),V(r,X(r)),V_x(r,X(r))^\top\si(r,X(r),u(r))\big)\\
\ns\ds\qq\qq\qq-V_t(r,X(r))-\dbH\big(r,X(r),u(r),V(r,X(r)),V_x(r,X(r)),V_{xx}(r,X(r))\big)\Big\}dr\\
\ns\ds\qq\qq\qq-\int_s^T V_x(r,X(r))^\top\si(r,X(r),u(r))dW(r),\qq s\in[t,T].\ea\ee
Noticing \eqref{Sec2.2_VT_Temp1}, by the Comparison Theorem of BSDEs
(see El Karoui--Peng--Quenez \cite{El Karoui-Peng-Quenez1997}), we
obtain
$$V(s,X(s))\les Y(s;t,\xi,u(\cd)),\qq\as,\q\forall\ s\in [t,T].$$
Particularly, by taking $s=t$, the above implies \eqref{Sec2.2_VT_1}.

\ms

Next, if the admissible pair $(\bar X(\cd),\bar u(\cd))$
satisfies \eqref{Sec2.2_VT_2}, then the equal sign in
\eqref{Sec2.2_VT_Temp1} holds, and \eqref{Sec2.2_VT_Temp2} becomes
$$\ba{ll}
\ns\ds V(s,\bar X(s))=h(\bar X(T))+\int_s^Tg\(r,\bar X(r),\bar
u(r),V(r,\bar X(r)),V_x(r,\bar X(r))^\top\si(r,\bar X(r),\bar u(r))\)dr\\
\ns\ds\qq\qq\qq\qq\qq-\int_s^T V_x(r,\bar X(r))^\top\si(r,\bar X(r),\bar u(r))dW(r)\qq s\in[t,T].\ea$$
By the uniqueness of adapted solutions to BSDE, we have
$$V(s,\bar X(s))=Y(s;t,\xi,\bar u(s)),\qq\as,\ \forall\ s\in [t,T].$$
Letting $s=t$, we obtain \eqref{Sec2.2_VT_3} and the proof is completed. \endpf

\section{Time-Inconsistent Problem and Equilibrium Strategy}\label{Sec_Formulation}

In this section, we formulate the time-inconsistent recursive stochastic optimal
control problem, and introduce the notion of equilibrium strategy for the problem.

\ms

For convenience, let us rewrite the state equation \eqref{Sec1_FBSDE_Inc} and the recursive cost functional \eqref{recursive} as follows:
\bel{Sec3_FBSDE}\left\{\2n\ba{ll}
\ns\ds dX(s)=b(s,X(s),u(s))ds+\si(s,X(s),u(s))dW(s),\qq s\in[t,T],\\
\ns\ds dY(s)=-g(t,s,X(s),u(s),Y(s),Z(s))ds+Z(s)dW(s),\qq s\in[t,T],\\
\ns\ds X(t)=\xi,\qq Y(T)=h(t,X(T)),\ea\right.\ee
\bel{Sec3_Cost}J(t,\xi;u(\cd))=Y(t;t,\xi,u(\cd)),\qq u(\cd)\in\sU[t,T].\ee
We pose the following problem.

\ms

{\bf Problem (N).} For any $(t,\xi)\in\sD$, find a $\bar u(\cd)\in\sU[t,T]$ such that
$$J(t,\xi;\bar u(\cd))=\essinf_{u(\cd)\in\sU[t,T]}J(t,\xi;u(\cd)).$$

As mentioned in Section \ref{Sec_Introduction}, Problem (N) is
time-inconsistent. Therefore, instead of find an optimal control for a given initial pair $(t,\xi)\in\sD$, we would like to find a time-consistent equilibrium strategy for Problem (N) over the whole time interval $[0,T]$. Our approach is inspired by that developed in \cite{Yong2012b}, for a time-inconsistent optimal control problem of stochastic differential equations (SDE, for short), with Bolza type cost functional (We now have a recursive cost functional).

\ms

Let $\cP[0,T]$ denote the set of all partitions $\Pi=\{t_k\
|\ 0\les k\les N\}$ of $[0,T]$ with $0=t_0<t_1<t_2<\cds<t_{N-1}<t_N=T$.
The mesh size of $\Pi$ is defined as $\|\Pi\|=\ds\max_{1\les k\les N}(t_k-t_{k-1})$. Similar to \cite{Yong2012b}, we present the definition of time-consistent equilibrium strategy of
Problem (N) as follows.

\begin{definition}\label{Sec3_Def} \rm A continuous map $\mathbbm u:[0,T]\times\dbR^n\to U$
is called a {\rm time-consistent equilibrium strategy} of Problem
(N) if the following hold:

\ms

(i) (Time-consistency) For any $x\in\dbR^n$, the following closed-loop system:
\bel{closed-loop}\left\{\2n\ba{ll}
\ns\ds d\bar X(s)=b\big(s,\bar X(s),\mathbbm u(s,\bar X(s))\big)ds
+\si\big(s,\bar X(s),\mathbbm u(s,\bar X(s))\big)dW(s),\qq s\in [0,T],\\
\ns\ds d\bar Y(t,s)=-g\big(t,s,\bar X(s),\mathbbm u(s,\bar
X(s)),\bar Y(t,s),\bar Z(t,s)\big)ds+Z(t,s)dW(s),\qq(t,s)\in D[0,T],\\
\ns\ds\bar X(0)=x,\qq\bar Y(t,T)=h\big(t,\bar X(T)\big),\qq t\in[0,T],\ea\right.\ee
admits a unique adapted solution $\big(\bar X(\cd),\bar Y(\cd\,,\cd),\bar Z(\cd\,,\cd)\big)\equiv
\big(\bar X(\cd\,;x),\bar Y(\cd\,,\cd\,;x),Z(\cd\,,\cd\,;x)\big)$.

\ms

(ii) (Local approximate optimality) There exists a family of partitions $\cP_0[0,T]\subseteq\cP[0,T]$ with
$$\inf_{\Pi\in\cP_0[0,T]}\|\Pi\|=0,$$
and a family of maps $\mathbbm u^\Pi:[0,T]\times\dbR^n\to U$ parameterized by the partitions $\Pi\in\cP_0[0,T]$ such that
$$\lim_{\|\Pi\|\to0}\sup_{(t,x)\in\cK}|\mathbbm u^\Pi(t,x)-\mathbbm u(t,x)|=0,\qq\forall\hb{ compact set
$\cK\subseteq[0,T]\times\dbR^n$},$$
with $\mathbbm u^\Pi$ being locally optimal in the following sense: Let $\Pi=\{t_k\
|\ 0\les k\les N\}\in\cP_0[0,T]$. For any $x\in\dbR^n$, the following FBSDE:
$$\left\{\2n\ba{ll}
\ns\ds dX^\Pi(s)=b\big(s,X^\Pi(s),\mathbbm u^\Pi(s,X^\Pi(s))\big)ds
+\si\big(s,X^\Pi(s),\mathbbm u^\Pi(s,X^\Pi(s))\big)dW(s),\qq s\in[0,T],\\
\ns\ds dY^\Pi(t_k,s)=-g\big(t_k,s,X^\Pi(s),\mathbbm u^\Pi(s,X^\Pi(s)),
Y^\Pi(t_k,s),Z^\Pi(t_k,s)\big)ds+Z^\Pi(t_k,s)dW(s),\\
\ns\ds\qq\qq\qq\qq\qq\qq\qq\qq\qq\qq s\in [t_k,T],\qq0\les k\les N-1,\\
\ns\ds X^\Pi(0)= x,\qq Y^\Pi(t_k,T)=h\big(t_k,X^\Pi(T)\big),\qq0\les k\les N-1\ea\right.$$
admits a unique adapted solution:
$$(X^\Pi(\cd),Y^\Pi(t_k,\cd),Z^\Pi(t_k,\cd))\equiv
(X^\Pi(\cd\,;x,\mathbbm u^\Pi(\cd)),Y^\Pi(t_k,\cd\,;x,\mathbbm u^\Pi(\cd)),Z^\Pi(t_k,\cd\,;x,\mathbbm u^\Pi(\cd))),$$
such that for each $k=1,2,\cds,N$,
\bel{Sec3_Nash}\ba{ll}
\ns\ds J\big(t_{k-1},X^\Pi(t_{k-1});\mathbbm
u^\Pi(\cd\,;X^\Pi(\cd))|_{[t_{k-1},T]}\big)=Y^\Pi(t_{k-1},t_{k-1})\\
\ns\ds\qq\les J\big(t_{k-1},X^\Pi(t_{k-1});u^k(\cd)\oplus\mathbbm
u^\Pi(\cd\,;X^k(\cd))|_{[t_k,T]}\big)= Y^k(t_{k-1}),\qq \forall\
u^k(\cd)\in\sU[t_{k-1},t_k],\ea\ee
where
\bel{Sec3_oplus}\[u^k(\cd)\oplus\mathbbm u^\Pi(\cd\,;X^k(\cd))|_{[t_k,T]}\](s)=
\left\{\2n\ba{ll}
\ns\ds u^k(s),\qq\qq\qq s\in[t_{k-1},t_k),\\
\ns\ds\mathbbm u^\Pi\big(s,X^k(s)\big),\qq\q s\in[t_k,T],\ea\right.\ee
and $(X^k(\cd),Y^k(\cd),Z^k(\cd))$ is the unique adapted solution to
the following FBSDE defined on the interval $[t_{k-1},T]$:
$$\left\{\2n\ba{ll}
\ns\ds dX^k(s)=b\big(s,X^k(s),u^k(s)\big)ds+\si\big(s,X^k(s),u^k(s)\big)dW(s),\qq s\in[t_{k-1},t_k),\\
\ns\ds dX^k(s)=b\big(s,X^k(s),\mathbbm u^\Pi(s,X^k(s))\big)ds+\si\big(s,X^k(s),\mathbbm u^\Pi(s,X^k(s)) \big)dW(s),\qq s\in [t_k,T],\\
\ns\ds dY^k(s)=-g\big(t_{k-1},s,X^k(s),u^k(s),Y^k(s),Z^k(s)\big)ds+Z^k(s)dW(s),\qq
s\in [t_{k-1},t_k),\\
\ns\ds dY^k(s)=-g\big(t_{k-1},s,X^k(s),\mathbbm u^\Pi(s,X^k(s)),Y^k(s),Z^k(s)\big)ds+Z^k(s)dW(s),\qq s\in [t_k,T],\\
\ns\ds X^k(t_{k-1})=X^\Pi(t_{k-1}),\qq Y^k(T)=h\big(t_{k-1},X^k(T)\big).\ea\right.$$
Moreover, for any $(s,x)\in[0,T]\times\dbR^n$,
\bel{Sec3-Convergence}\lim_{\|\Pi\|\to0}\(|X^\Pi(s)-\bar X(s)|+|Y^\Pi(\ell^\Pi(s),\ell^\Pi(s))-\bar Y(s,s)|\)=0,\qq\as,\ee
where
\bel{Sec3_l_Pi}\ell^\Pi(s)=\sum_{k=1}^Nt_{k-1}\mathbbm 1_{[t_{k-1},t_k)}(s),\qq
s\in [0,T].\ee

In the above, $\bar X(\cd)$ is called a {\it time-consistent
equilibrium state process}, $\bar u(\cd)\equiv\mathbbm
u(\cd\,,\bar X(\cd))$ is called a {\it time-consistent equilibrium control}
for the initial state $x$, and $(\bar X(\cd),\bar u(\cd))$ is called
a {\it time-consistent equilibrium pair} of Problem (N). It is easy to see
that the convergence \eqref{Sec3-Convergence} implies
$$\lim_{\|\Pi\|\to0}J\big(\ell^\Pi(t),X^\Pi(\ell^\Pi(t));\mathbbm
u^\Pi(\cd\,,X^\Pi(\cd))|_{[\ell^\Pi(t),T]}\big)=J\big(t,\bar X(t),\mathbbm
u(\cd\,,\bar X(\cd))|_{[t,T]}\big),\q\as$$
We call $V:[0,T]\times\dbR^n\to\dbR^n$ an {\it equilibrium value function} of Problem (N) if
\bel{Sec3_ValueFucntion}V\big(t,\bar X(t)\big)=J\big(t,\bar X(t),\mathbbm
u(\cd\,,\bar X(\cd))|_{[t,T]}\big),\q\as, \ (t,x)\in
[0,T]\times\dbR^n.\ee
We also call $\mathbbm u^\Pi(\cd\,,\cd)$ an {\it approximate equilibrium strategy} of Problem (N) associated with the partition $\Pi$.

\end{definition}

Let us make a couple of comments on the above long-looking definition.

\ms

$\bullet$ The state equation \eqref{closed-loop} admits a unique solution under $\mathbbm u(\cd\,,\cd)$ means that as a strategy, $\mathbbm u(\cd\,,\cd)$ is time-consistent. It is interesting to know that such a strategy is of closed-loop nature, in the sense that it is independent of the initial state.

\ms

$\bullet$ Condition \eqref{Sec3_Nash} means that the {\it outcome} $\mathbbm u^\Pi(\cd\,;X^\Pi(\cd))$ of the strategy $\mathbbm u^\Pi(\cd\,,\cd)$ is locally optimal in a proper sense. Due to the fact that the global optimal control is time-inconsistent, such kind of local optimal control should be the best that one can obtain.

\section{Multi-Person Differential Games}\label{Sec_Game}

In this section, we shall construct a family $\mathbbm u^\Pi(\cd\,,\cd)$ of approximate equilibrium strategies which plays the role as in Definition \ref{Sec3_Def}. To this end, we consider a family of multi-person differential games, called Problem (G$^\Pi$), associated with the partition $\Pi:0=t_0<t_1<\cds<t_{N-1}<t_N=T$. In the game, there are $N$ players labeled from $1$ to $N$. Player $k$ controls the system on the interval $[t_{k-1},t_k)$ by selecting his/her own admissible control $u^k(\cd)\in\sU[t_{k-1},t_k]$. We now carry out the details below.

\subsection{Player $N$ --- a classical optimal control problem}

In what follows, we denote
$$L^2_{\cF_t}(\O;\dbR^n)=\Big\{\xi:\O\to\dbR^n\bigm|\xi\hb{ is $\cF_t$-measurable, }\dbE|\xi|^2<\infty\Big\},\qq t\in[0,T].$$
Let us start with Player $N$ who controls the system on $[t_{N-1},t_N]$, the last time interval of the partition $\Pi$. For any admissible control $u^N(\cdot)\in\mathscr U[t_{N-1},t_N]$ and initial state $\xi_{N-1}\in L^2_{\cF_{t_{N-1}}}(\O;\dbR^n)$, the controlled FBSDE for Player $N$ reads
\bel{Sec4.1_FBSDE}\left\{\2n\ba{ll}
\ns\ds dX^N(s)=b\big(s,X^N(s),u^N(s)\big)ds +\si\big(s,X^N(s),u^N(s)\big)dW(s),\qq s\in[t_{N-1},t_N],\\
\ns\ds dY^N(s)=-g\big(t_{N-1},s,X^N(s),u^N(s),Y^N(s),Z^N(s)\big)ds+Z^N(s)dW(s),\qq s\in[t_{N-1},t_{N}],\\
\ns\ds X^N(t_{N-1})=\xi_{N-1},\qq Y^N(t_N)=h\big(t_{N-1},X^N(t_N)\big),\ea\right.\ee
whose unique adapted solution is denoted by
$$(X^N(\cd),Y^N(\cd),Z^N(\cd))\equiv\big(X^N(\cd\,;\xi_{N-1},u^N(\cd)),Y^N(\cd\,;\xi_{N-1},u^N(\cd)),
Z^N(\cd\,;\xi_{N-1},u^N(\cd))\big),$$
emphasizing the dependence on $(\xi_{N-1},u^N(\cd))$. The recursive cost functional is given by
\bel{Sec4.1_Cost}\ba{ll}
\ns\ds J(t_{N-1},\xi_{N-1};u^N(\cd))=Y^N(t_{N-1};\xi_{N-1},u^N(\cd))\\
\ns\ds\qq\qq=\dbE_{t_{N-1}}\[\int_{t_{N-1}}^Tg(t_{N-1},s,X^N(s),u^N(s),Y^N(s),Z^N(s))ds+h(t_{N-1},X^N(T))\].
\ea\ee
The optimal control problem for Player $N$ can be stated as follows.

\ms

{\bf Problem (C$_N$).} For any $\xi_{N-1}\in L^2_{\cF_{t_{N-1}}}(\O;\dbR^n)$, find a $\bar u^N(\cd)\in\sU[t_{N-1},t_N]$ such that
$$J(t_{N-1},\xi_{N-1};\bar u^N(\cd))=\essinf_{u^N(\cd)\in\sU[t_{N-1},t_N]}
J(t_{N-1},\xi_{N-1};u^N(\cd)).$$

Problem (C$_N$) is a standard stochastic optimal control problem with recursive cost functional, which can be solved by the stochastic verification theorem (see Theorem \ref{Sec2.2_THM_Verification}). More precisely, under proper conditions, the following HJB equation admits a classical solution $V^\Pi(\cd\,,\cd)$:
\bel{Sec4.1_HJB_inf}\left\{\2n\ba{ll}
\ns\ds V^\Pi_t(t,x)+\inf_{u\in U}\dbH\big(t_{N-1},t,x,u,V^\Pi(t,x),V^\Pi_x(t,x),V^\Pi_{xx}(t,x)\big)
=0,\qq(t,x)\in[t_{N-1},t_N]\times\dbR^n,\\
\ns\ds V^\Pi(t_N,x)=h(t_{N-1},x),\ea\right.\ee
where $\dbH$ is defined by \eqref{Sec2_Notations}. Recalling the map $\psi$ introduced in (H2), we define
\bel{Sec4.1_mathbbm_u}\mathbbm u^\Pi(t,x)=\psi\big(t_{N-1},t,x,V^\Pi(t,x),V^\Pi_x(t,x),V^\Pi_{xx}(t,x)\big),\qq
(t,x)\in [t_{N-1},t_N]\times\dbR^n.\ee
Let us assume that for $\xi_{N-1}\in L^2_{\cF_{t_{N-1}}}(\O;\dbR^n)$, the following FBSDE (which is a closed-loop system):
\bel{}\left\{\2n\ba{ll}
\ns\ds d\bar X^N(s)=b\big(s,\bar X^N(s),\mathbbm u^\Pi(s,\bar X^N(s))\big)ds+\si\big(s,\bar X^N(s),\mathbbm u^\Pi(s,\bar X^N(s))\big)dW(s),\q s\in [t_{N-1},t_N],\\
\ns\ds d\bar Y^N(s)=-g\big(t_{N-1},s,\bar X^N(s),\mathbbm u^\Pi(s,\bar X^N(s)),\bar Y^N(s),\bar Z^N(s)\big)ds+\bar Z^N(s)dW(s),\q s\in [t_{N-1},t_{N}],\\
\ns\ds\bar X^N(t_{N-1})=\xi_{N-1},\qq\bar Y^N(t_N)=h\big(t_{N-1},\bar X^N(t_N)\big)\ea\right.\ee
admits a unique adapted solution
$$\big(\bar X^N(\cd),\bar Y^N(\cd),\bar Z^N(\cd)\big)\equiv\big(\bar X^N(\cd\,;\xi_{N-1}),\bar Y^N(\cd\,;\xi_{N-1}),\bar Z^N(\cd\,;\xi_{N-1})\big).$$
Then by Theorem \ref{Sec2.2_THM_Verification}, $(\bar X^N(\cd\,;\xi_{N-1}),\mathbbm
u^\Pi(\cd\,,\bar X^N(\cd\,;\xi_{N-1})))$ is an optimal pair of Problem (C$_N$) for the initial pair $(t_{N-1},\xi_{N-1})$. Because of that, $\mathbbm u^\Pi(\cd\,,\cd)$ defined by \eqref{Sec4.1_mathbbm_u}
(on $[t_{N-1},t_N]\times\dbR^n$) is called an optimal strategy of Player $N$, and $\mathbbm u^\Pi(\cd\,,\bar X^N(\cd))$ is called an {\it outcome} of $\mathbbm u^\Pi(\cd\,,\cd)$.

\subsection{Player ($N-1$) --- a sophisticated optimal control problem}

We now look at Player ($N-1$) who takes over the system from Player $(N-2)$, controls the system on
$[t_{N-2},t_{N-1})$ and hand it over to Player $N$ at $t_{N-1}$. Player ($N-1$) knows that
Player $N$ will play optimally through the optimal strategy $\mathbbm
u^\Pi(\cd\,,\cd)$ (which is already defined on $[t_{N-1},t_N]\times\dbR^n$). Due to the subjective time-preference, Player ($N-1$) still
``discounts" the future costs in his/her own way despite he/she will not control the system beyond $t_{N-1}$. According to this viewpoint, the controlled FBSDE of Player $(N-1)$ is given by
\bel{Sec4.2_FBSDE_Multi}\left\{\2n\ba{ll}
\ns\ds dX^{N-1}(s)=b\big(s,X^{N-1}(s),u^{N-1}(s)\big)ds
+\si\big(s,X^{N-1}(s),u^{N-1}(s)\big)dW(s),\qq s\in
[t_{N-2},t_{N-1}),\\
\ns\ds dX^{N-1}(s)=b\big(s,X^{N-1}(s),\mathbbm
u^\Pi(s,X^{N-1}(s))\big)ds+\si\big(s,X^{N-1}(s),\mathbbm u^\Pi(s,X^{N-1}(s))\big)dW(s),\\
\ns\ds\hskip 11.8cm s\in [t_{N-1},t_N],\\
\ns\ds dY^{N-1}(s)=-g\big(t_{N-2},s,X^{N-1}(s),u^{N-1}(s),Y^{N-1}(s),Z^{N-1}(s)\big)ds
+Z^{N-1}(s)dW(s),\\
\ns\ds\hskip 11.8cm s\in [t_{N-2},t_{N-1}),\\
\ns\ds dY^{N-1}(s)=-g\big(t_{N-2},s,X^{N-1}(s),\mathbbm
u^\Pi(s,X^{N-1}(s)),Y^{N-1}(s),Z^{N-1}(s)\big)ds+Z^{N-1}(s)dW(s),\\
\ns\ds\hskip 11.8cm s\in [t_{N-1},t_N],\\
\ns\ds X^{N-1}(t_{N-2})=\xi_{N-2},\qq Y^{N-1}(t_N)=h\big(t_{N-2},X^{N-1}(t_N)\big),\ea\right.\ee
where $\xi_{N-2}\in L^2_{\cF_{t_{N-2}}}(\O;\dbR^n)$ and $u^{N-1}(\cd)\in\sU[t_{N-2},t_{N-1}]$.
Let
$$\ba{ll}
\ns\ds(X^{N-1}(\cd),Y^{N-1}(\cd),Z^{N-1}(\cd))\\
\ns\ds\equiv\big(X^{N-1}(\cd\,;\xi_{N-2},u^{N-1}(\cd)),Y^{N-1}(\cd\,;\xi_{N-2},u^{N-1}(\cd)),
Z^{N-1}(\cd\,;\xi_{N-2},u^{N-1}(\cd))\big)\ea$$
be the adapted solution of \eqref{Sec4.2_FBSDE_Multi}, depending on the initial state $\xi_{N-2}\in L^2_{\cF_{t_{N-2}}}(\O;\dbR^n)$ and the control $u^{N-1}(\cd)\in\sU[t_{N-2},t_{N-1}]$. Then we define the {\it sophisticated} recursive cost functional of Player $(N-1)$ by the following:
\bel{Sec4.2cost}\wt J(t_{N-2},\xi_{N-2};u^{N-1}(\cd))\equiv J(t_{N-2},\xi_{N-2};u^{N-1}(\cd)\oplus\mathbbm
u^\Pi(\cd\,,X^{N-1}(\cd)))=Y^{N-1}(t_{N-2};\xi_{N-2},u^{N-1}(\cd)).\ee
where the operation ``$\oplus$" is defined by \eqref{Sec3_oplus} and $\mathbbm u^\Pi(\cd\,,\cd)$ is defined by \eqref{Sec4.1_mathbbm_u}. Clearly, the above is different from the ``naive'' recursive cost functional $J(t_{N-2},\xi_{N-2};u(\cd))$ defined by \eqref{recursive}. We emphasize that for different control $u^{N-1}(\cd)$ selected from $\sU[t_{N-2},t_{N-1}]$, $X^{N-1}(t_{N-1})\equiv X^{N-1}(t_{N-1};t_{N-2},\xi_{N-2},u^{N-1}(\cd))$ will be different, which will result in the FBSDE on $[t_{N-1},t_N]$ having a different initial condition for $X^{N-1}(\cd)$. Now we pose the following problem for Player $(N-1)$.

\ms

{\bf Problem (C$_{N-1}$).} For any $\xi_{N-2}\in L^2_{\cF_{t_{N-2}}}(\O;\dbR^n)$, find a $\bar
u^{N-1}(\cd)\in\sU[t_{N-2},t_{N-1}]$ such that
\bel{J(N-1)}\wt J\big(t_{N-2},\xi_{N-2};\bar u^{N-1}(\cd)\big)=\essinf_{u^{N-1}(\cd)\in\sU[t_{N-2},t_{N-1}]}
\wt J\big(t_{N-2},\xi_{N-2},u^{N-1}(\cd)\big).\ee

\ms

Let us make some careful observation which will reveal the essential difference between the naive and sophisticated recursive cost functionals. Recall that the original controlled FBSDE on $[t_{N-2},t_N]$ is given by
\bel{Sec4.2_FBSDE}\left\{\2n\ba{ll}
\ns\ds dX(s)=b(s,X(s),u(s))ds+\si(s,X(s),u(s))dW(s),\qq s\in[t_{N-2},T],\\
\ns\ds dY(s)=-g(t_{N-2},s,X(s),u(s),Y(s),Z(s))ds+Z(s)dW(s),\qq s\in[t_{N-2},T],\\
\ns\ds X(t_{N-2})=\xi_{N-2},\qq Y(T)=h(t_{N-2},X(T)),\ea\right.\ee
and the naive recursive cost functional is given by
\bel{Sec4_Cost}J(t_{N-2},\xi_{N-2};u(\cd))=Y(t_{N-2};t_{N-2},\xi,u(\cd)),\qq u(\cd)\in\sU[t_{N-2},T].\ee
If $t_{N-2}$ is fixed (as a parameter), the above will lead to a time-consistent optimal control problem with recursive cost functional. Consequently, if $(\wt X^{N-2}(\cd),\wt u^{N-2}(\cd))$ is an optimal pair (corresponding to the initial pair $(t_{N-2},\xi_{N-2})$), then
\bel{wt u}\ba{ll}
\ns\ds\wt u^{N-2}(s)=\wt{\mathbbm u}^{N-2}(s,\wt X^{N-2}(s)),\qq s\in[t_{N-2},T],\\
\ns\ds\wt{\mathbbm u}^{N-2}(t,x)=\psi\big(t_{N-2},t,x,\wt V^{N-2}(t,x),\wt V^{N-2}_x(t,x),\wt V^{N-2}_{xx}(t,x)\big),\qq(t,x)\in[t_{N-2},T]\times\dbR^n,\ea\ee
with $\wt V^{N-2}(\cd\,,\cd)$ satisfying the following HJB equation:
\bel{Sec4.2_HJB_inf}\left\{\2n\ba{ll}
\ns\ds\wt V^{N-2}_t(t,x)+\inf_{u\in U}\dbH\big(t_{N-2},t,x,u,\wt V^{N-2}(t,x),\wt V^{N-2}_x(t,x),\wt V^{N-2}_{xx}(t,x)\big)=0,\\
\ns\ds\qq\qq\qq\qq\qq\qq\qq\qq\qq\qq\qq\qq\qq(t,x)\in[t_{N-2},T]\times\dbR^n,\\
\ns\ds\wt V^{N-2}(t_N,x)=h(t_{N-2},x),\qq x\in\dbR^n.\ea\right.\ee
Comparing \eqref{Sec4.1_HJB_inf} and \eqref{Sec4.2_HJB_inf}, we see that on $[t_{N-1},T]\times\dbR^n$, $V^\Pi(\cd\,,\cd)$ and $\wt V^{N-2}(\cd\,,\cd)$ satisfy different HJB equations: The former has the parameter $t_{N-1}$ and the later has the parameter $t_{N-2}$. Hence, they are different in general. Consequently, by further comparing \eqref{Sec4.1_mathbbm_u} and \eqref{wt u}, we see that $\mathbbm u^\Pi(\cd\,,\cd)$ and $\wt{\mathbbm u}^{N-2}(\cd\,,\cd)$ are different. In another word,
$$\wt u^{N-2}(\cd)=\bar u^{N-1}(\cd)\oplus\mathbbm u^\Pi(\cd\,;\bar X^{N-1}(\cd))$$
fails in general. The right hand side of the above is called an equilibrium control of Problem (G$^\Pi$) on $[t_{N-2},T]$, which is not an optimal control of Problem (N) on $[t_{N-2},T]$.

\ms

Now, we would like to obtain a better representation of the sophisticated cost functional \eqref{Sec4.2cost} of Player $(N-1)$. To this end, we look at the following closed-loop system on $[t_{N-1},t_N]$:
\bel{Sec4.2_FBSDE_B}\left\{\2n\ba{ll}
\ns\ds dX^{N-1}(s)=b\big(s,X^{N-1}(s),\mathbbm u^\Pi(s,X^{N-1}(s))\big)ds
+\si\big(s,X^{N-1}(s),\mathbbm u^\Pi(s,X^{N-1}(s))\big)dW(s),\\
\ns\ds\hskip 11.8cm s\in [t_{N-1},t_N],\\
\ns\ds dY^{N-1}(s)=-g\big(t_{N-2},s,X^{N-1}(s),\mathbbm
u^\Pi(s,X^{N-1}(s)),Y^{N-1}(s),Z^{N-1}(s)\big)ds+Z^{N-1}(s)dW(s),\\
\ns\ds\hskip 11.8cm s\in [t_{N-1},t_N],\\
\ns\ds X^{N-1}(t_{N-1})=X^{N-1}(t_{N-1}),\qq Y^{N-1}(t_N)=h\big(t_{N-2},X^{N-1}(t_N)\big).\ea\right.\ee
Inspired by Theorem \ref{Sec2.1_THM_Feynman-Kac}, we introduce the following
PDE:
\bel{Sec4.2_PDE}\left\{\2n\ba{ll}
\ns\ds\Th^{N-1}_t(t,x)+\mathbb H\big(t_{N-2},t,x,\mathbbm
u^\Pi(t,x),\Th^{N-1}(t,x),\Th^{N-1}_x(t,x),\Th^{N-1}_{xx}(t,x)\big)=0,\\
\ns\ds\hskip 8cm (t,x)\in [t_{N-1},t_N]\times\dbR^n,\\
\ns\ds\Th^{N-1}(t_N,x)=h(t_{N-2},x),\qq x\in\dbR^n.\ea\right.\ee
If the above admits a unique classical solution $\Th^{N-1}(\cd\,,\cd)$, then the following representation holds:
\bel{Y(N-1)=Th(N-1)}Y^{N-1}(s)=\Th^{N-1}(s,X^{N-1}(s)),\qq s\in[t_{N-1},t_N].\ee
In particular,
\bel{Y(N-1)=Th(N-1)*}Y^{N-1}(t_{N-1})=\Th^{N-1}\big(t_{N-1},X^{N-1}(t_{N-1})\big).\ee
Consequently, the sophisticated cost functional will have the representation
\bel{wt J(N-2)}\wt J(t_{N-2},\xi_{N-2};u^{N-1}(\cd))=Y^{N-1}(t_{N-2};\xi_{N-2},u^{N-1}(\cd)),\ee
with
$$\ba{ll}
\ns\ds\big(X^{N-1}(\cd\,;\xi_{N-2},u^{N-1}(\cd)),Y^{N-1}(\cd\,;\xi_{N-2},u^{N-1}(\cd)),
Z^{N-1}(\cd\,;\xi_{N-2},u^{N-1}(\cd))\big)\\
\ns\ds\equiv(X^{N-1}(\cd),Y^{N-1}(\cd),Z^{N-1}(\cd))\ea$$
being the adapted solution to the following decoupled FBSDE:
\bel{}\left\{\2n\ba{ll}
\ns\ds dX^{N-1}(s)=b\big(s,X^{N-1}(s),u^{N-1}(s)\big)ds
+\si\big(s,X^{N-1}(s),u^{N-1}(s)\big)dW(s),\qq s\in
[t_{N-2},t_{N-1}),\\
\ns\ds dY^{N-1}(s)=-g\big(t_{N-2},s,X^{N-1}(s),u^{N-1}(s),Y^{N-1}(s),Z^{N-1}(s)\big)ds
+Z^{N-1}(s)dW(s),\\
\ns\ds\hskip 11.8cm s\in [t_{N-2},t_{N-1}),\\
\ns\ds X^{N-1}(t_{N-2})=\xi_{N-2},\qq Y^{N-1}(t_{N-1})=\Th^{N-1}\big(t_{N-1},X^{N-1}(t_{N-1})\big).\ea\right.\ee
Then Problem (C$_{N-1}$) becomes a standard recursive stochastic optimal control problem (on $[t_{N-2},t_{N-1}]$). Let the following HJB equation admits a classical solution:
\bel{}\left\{\2n\ba{ll}
\ns\ds V^\Pi_t(t,x)+\inf_{u\in U}\dbH\big(t_{N-2},t,x,u,V^\Pi(t,x),V^\Pi_x(t,x),V^\Pi_{xx}(t,x)\big)=0,\q(t,x)\in [t_{N-2},t_{N-1})\times\mathbb R^n,\\
\ns\ds V^\Pi(t_{N-1}-0,x)=\Th^{N-1}(t_{N-1},x),\qq x\in\dbR^n.\ea\right.\ee
Similar to \eqref{Sec4.1_mathbbm_u}, we define
\bel{Sec4.2_mathbbm_u_F}\mathbbm u^\Pi(t,x)=
\psi\big(t_{N-2},t,x,V^\Pi(t,x),V^\Pi_x(t,x),V^\Pi_{xx}(t,x)\big),\qq
(t,x)\in [t_{N-2},t_{N-1})\times\dbR^n.\ee
We assume the following FBSDE
\bel{}\left\{\2n\ba{ll}
\ns\ds d\bar X^{N-1}(s)=b\big(s,\bar X^{N-1}(s),u^\Pi(s,\bar X^{N-1}(s))\big)ds
+\si\big(s,\bar X^{N-1}(s),u^\Pi(s,\bar X^{N-1}(s))\big) dW(s),\\
\ns\ds\hskip 11cm s\in[t_{N-2},t_{N-1}],\\
\ns\ds d\bar Y^{N-1}(s)=-g\big( t_{N-2},s,\bar X^{N-1}(s),u^\Pi(s,\bar X^{N-1}(s)),\bar Y^{N-1}(s),\bar Z^{N-1}(s)\big)ds+\bar Z^{N-1}(s)dW(s),\\
\ns\ds\hskip 11cm s\in [t_{N-2},t_{N-1}],\\
\ns\ds\bar X^{N-1}(t_{N-2})=\xi_{N-2},\qq\bar Y^{N-1}(t_{N-1})=\Th^{N-1}\big(t_{N-1},\bar
X^{N-1}(t_{N-1})\big)\ea\right.\ee
admits a unique adapted solution. Then by Theorem \ref{Sec2.2_THM_Verification},
$(\bar X^{N-1}(\cd),\mathbbm u^\Pi(\cd\,,\bar X^{N-1}(\cd)))$ is an optimal pair of Problem
(C$_{N-1}$). The map $\mathbbm u^\Pi(\cd\,,\cd)$ defined by \eqref{Sec4.2_mathbbm_u_F} is called an optimal strategy of Player ($N-1$).

\ms

Now, combining the optimal strategies and value functions of Players $N$ and $(N-1)$, we obtain that
both $\mathbbm u^\Pi(\cd\,,\cd)$ and $V^\Pi(\cd\,,\cd)$ are defined on $[t_{N-2},t_N]\times\dbR^n$, with a possible jump at $t=t_{N-1}$.

We write \eqref{Sec4.1_mathbbm_u} and \eqref{Sec4.2_mathbbm_u_F}
compactly as
\bel{Sec4.2_mathbbm_u}\mathbbm u^\Pi(t,x)=\psi\big(\ell^\Pi(t),t,x,V^\Pi(t,x),V^\Pi_x(t,x),V^\Pi_{xx}(t,x)\big),\q(t,x)\in \([t_{N-2},t_{N-1})\cup(t_{N-1},t_N]\)\times\dbR^n,\ee
where $\ell^\Pi(\cd)$ is defined by \eqref{Sec3_l_Pi}. From the above, we see that, respectively restricted in $[t_{N-2},t_{N-1})$ and $(t_{N-1},t_N]$, $\mathbbm u^\Pi(\cd\,,\cd)$ is an optimal strategy of Player ($N-1$) and Player $N$, respectively. However, in general $\mathbbm u^\Pi(\cd\,,\cd)$ is not an optimal strategy on the whole interval $[t_{N-2},t_N]$. We call $\mathbbm u^\Pi(\cd\,,\cd)$ an {\it equilibrium strategy} of Problem (G$^\Pi$) on $[t_{N-2},t_N]$. Further,
$$\big(\bar X^{N-1}(s),\mathbbm u^\Pi(s,\bar X^{N-1}(s))\big),\qq s\in [t_{N-2},t_N]$$
is called an {\it equilibrium pair} on $[t_{N-2},t_N]$, where
$$\bar X^{N-1}(s)=\left\{\2n\ba{ll}
\ns\ds\bar X^N(s),\qq\q s\in [t_{N-1},t_N],\\
\ns\ds\bar X^{N-1}(s),\qq s\in [t_{N-2},t_{N-1}).\ea\right.$$

\subsection{Player $k$ and equilibria of Problem (G$^\Pi$)}\label{SubSec_Player_k}

The above procedure can be continued recursively. Suppose we have constructed the equilibrium strategy $\mathbbm u^\Pi(\cd\,,\cd)$ and the equilibrium value function $V^\Pi(\cd\,,\cd)$ on $[t_k,t_N]\times\dbR^n$ for Problem (G$^\Pi$). We now extend them to $[t_{k-1},t_k)\times\dbR^n$. On $[t_{k-1},t_k)$, Player $k$ controls the system and he/she knows that later players will play through the equilibrium strategy $\mathbbm u^\Pi(\cd\,,\cd)$, and meanwhile Player $k$ ``discounts" the
future costs in his/her own way. Hence, for any $\xi_{k-1}\in L^2_{\cF_{t_{k-1}}}(\O;\dbR^n)$ and $u^k(\cd)\in\sU[t_{k-1},t_k]$, the controlled FBSDE for Player $k$ reads:
\bel{Sec4.3-FBSDE(k)}\left\{\2n\ba{ll}
\ns\ds dX^k(s)=b\big(s,X^k(s),u^k(s)\big)ds+\si\big(s,X^k(s),u^k(s)\big)dW(s),\q s\in[t_{k-1},t_k),\\
\ns\ds dX^k(s)=b\big(s,X^k(s),\mathbbm u^\Pi(s,X^k(s))\big)ds
+\si\big(s,X^k(s),\mathbbm u^\Pi(s,X^k(s))\big)dW(s),\q s\in[t_k,t_N],\\
\ns\ds dY^k(s)=-g\big(t_{k-1},s,X^k(s),u^k(s),Y^k(s),Z^k(s)\big)ds+Z^k(s)dW(s),\q s\in[t_{k-1},t_k),\\
\ns\ds dY^k(s)=-g\big(t_{k-1},s,X^k(s),\mathbbm
u^\Pi(s,X^k(s)),Y^k(s),Z^k(s)\big)ds+Z^k(s)dW(s),\q s\in[t_k,t_N],\\
\ns\ds X^k(t_{k-1})=\xi_{k-1},\qq Y^k(t_N)=h\big(t_{k-1},X^k(t_N)\big).\ea\right.\ee
Suppose $(X^k(\cd),Y^k(\cd),Z^k(\cd))$ is the adapted solution to the above. Then the sophisticated recursive cost functional of Player $k$ is defined as follows:
$$\wt J(t_{k-1},\xi_{k-1};u^k(\cd))=J\big(t_{k-1},\xi_{k-1};u^k(\cd)\oplus\mathbbm u^\Pi(\cd\,,X^k(\cd))\big)=Y^k(t_{k-1}).$$
The optimal control problem of Player $k$ is given by

\ms

\noindent{\bf Problem (C$_k$).} For any $\xi_{k-1}\in L^2_{\cF_{t_{k-1}}}(\O:\dbR^n)$, find a $\bar
u^k(\cd)\in\sU[t_{k-1},t_k]$ such that
\bel{}\wt J(t_{k-1},\xi_{k-1};\bar u^k(\cd))=\essinf_{u^k(\cd)\in\sU[t_{k-1},t_k]}
\wt J(t_{k-1},\xi_{k-1};u^k(\cd)).\ee

To get a better representation for the sophisticated recursive cost functional of Player $k$, we introduce the following PDE associated with a part of \eqref{Sec4.3-FBSDE(k)}:
\bel{Sec4.3_PDE}\left\{\2n\ba{ll}
\ns\ds\Th^k_t(t,x)+\dbH\big(t_{k-1},t,x,\mathbbm
u^\Pi(t,x),\Th^k(t,x),\Th^k_x(t,x),\Th^k_{xx}(t,x)
\big)=0,\qq(t,x)\in[t_k,t_N]\times\dbR^n,\\
\ns\ds\Th^k(t_N,x)=h(t_{k-1},x),\qq x\in\dbR^n.\ea\right.\ee
Under proper conditions, PDE \eqref{Sec4.3_PDE} admits a classical
solution $\Th^k(\cd\,,\cd)$, and Theorem \ref{Sec2.1_THM_Feynman-Kac} leads to
\begin{equation}\label{Sec4.3_Y_Theta}
Y^{k}(t_{k}) = \Theta^{k}\big(t_{k},X^{k}(t_{k})\big).
\end{equation}

\begin{remark} \rm
From the definition, $\mathbbm u^\Pi(\cdot,\cdot)$ (see
\eqref{Sec4.2_mathbbm_u} for the case of $k=N-1$) may have jumps at $t_{k+1},t_{k+2},\dots,t_{N-1}$ in general. Due to this, by saying $\Th^k(\cd\,,\cd)$ being the classical solution of \eqref{Sec4.3_PDE}, we mean that $\Th^k(\cd\,,\cd)$ is continuous, on all intervals $[t_{N-1},t_N]$, $[t_{N-2},t_{N-1})$,...,$[t_k,t_{k+1})$, $\Th^k(\cd\,,\cd)$ is the classical solution, and $\Th^k_x(\cd\,,\cd)$ and $\Th^k_{xx}(\cd\,,\cd)$ are allowed to have jumps at $t_{N-1},t_{N-2},\cds,t_{k+1}$. For simplicity of notations, we write in a compact form \eqref{Sec4.3_PDE}.
\end{remark}

The FBSDE controlled by Player $k$ on the interval $[t_{k-1},t_{k}]$ is given by
\bel{}\left\{\2n\ba{ll}
\ns\ds dX^k(s)=b\big(s,X^k(s),u^k(s)\big)ds+\si\big(s,X^k(s),u^k(s)\big)dW(s),\qq s\in[t_{k-1},t_k],\\
\ns\ds dY^k(s)=-g\big(t_{k-1},s,X^k(s),u^k(s),Y^k(s),Z^k(s)\big)ds+Z^k(s)dW(s),\qq s\in[t_{k-1},t_k],\\
\ns\ds X^k(t_{k-1})=\xi_{k-1},\qq Y^k(t_k)=\Th^k\big(t_k,X^k(t_k)\big).\ea\right.\ee
Under proper conditions, the following HJB equation
\bel{Sec4.3_HJB}\left\{\2n\ba{ll}
\ns\ds V^\Pi_t(t,x)+\inf_{u\in U}\dbH\big(t_{k-1},t,x,u,V^\Pi(t,x),V^\Pi_x(t,x),V^\Pi_{xx}(t,x)\big)=0,\qq(t,x)\in[t_{k-1},t_k)\times\dbR^n,\\
\ns\ds V^\Pi(t_{k}-0,x)=\Th^k(t_k,x),\qq x\in\dbR^n,\ea\right.\ee
admits a classical solution $V^\Pi(\cd\,,\cd)\in C^{1,2}([t_{k-1},t_k)\times\dbR^n)$. Define
\bel{Sec4.3_mathbbm_u_F}\ba{ll}
\ns\ds\mathbbm u^\Pi(t,x)=\psi(t_{k-1},t,x,V^\Pi(t,x),V^\Pi_x(t,x),V^\Pi_{xx}(t,x)),\\
\ns\ds\qq\qq=\psi(\ell^\Pi(t),t,x,V^\Pi(t,x),V^\Pi_x(t,x),V^\Pi_{xx}(t,x)),
\qq(t,x)\in[t_{k-1},t_k)\times\dbR^n.\ea\ee
Moreover, we assume the following FBSDE
\bel{Sec4.3_FBSDE_F_Optimal}\left\{\2n\ba{ll}
\ns\ds d\bar X^k(s)=b\big(s,\bar X^k(s),\mathbbm u^\Pi(s,\bar X^k(s))
\big)ds+\si\big(s,\bar X^k(s),\mathbbm u^\Pi(s,\bar X^k(s))\big) dW(s),\qq s\in[t_{k-1},t_k],\\
\ns\ds d\bar Y^k(s)=-g\big(t_{k-1},s,\bar X^k(s),\mathbbm
u^\Pi(s,\bar X^k(s)),\bar Y^k(s),\bar Z^k(s)\big)ds+\bar Z^k(s)
dW(s),\qq s\in[t_{k-1},t_k],\\
\ns\ds\bar X^k(t_{k-1})=\xi_{k-1},\qq\bar Y^k(t_k)=\Th^k\big(t_k,\bar X^k(t_k)\big)\ea\right.\ee
admits a unique adapted solution. Then, by Theorem \ref{Sec2.2_THM_Verification}, $\mathbbm u^\Pi(\cd\,,\cd)$ defined by \eqref{Sec4.3_mathbbm_u_F} is an optimal strategy of Player $k$.
Further, both $V^\Pi(\cd\,,\cd)$ and $\mathbbm u^\Pi(\cd\,,\cd)$ are now defined on $[t_{k-1},t_N]\times\dbR^n$.

\ms

By induction, we are able to obtain $V^\Pi(\cd\,,\cd)$ and $\mathbbm u^\Pi(\cd\,,\cd)$ defined on $[0,t_N]\times\dbR^n$. We point out that the construction of $V^\Pi(\cd\,,\cd)$ and $\mathbbm u^\Pi(\cd\,,\cd)$ is recursive. For the later purpose of taking the limits, let us summarize the procedure of construction.

\ms

\it Step 1. \rm Define value function $V^\Pi(\cd\,,\cd)$ on $[t_{N-1},t_N]\times\dbR^n$ through HJB equation:
\bel{HJB(N)}\left\{\2n\ba{ll}
\ns\ds V^\Pi_t(t,x)+\inf_{u\in U}\dbH\big(t_{N-1},t,x,u,V^\Pi(t,x),V_x^\Pi(t,x),V^\Pi_{xx}(t,x)\big)=0,\qq
(t,x)\in[t_{N-1},t_N]\times\dbR^n,\\
\ns\ds V^\Pi(t_N,x)=h(t_{N-1},x),\qq x\in\dbR^n.\ea\right.\ee
Having $V^\Pi(\cd\,,\cd)$ defined on $[t_{N-1},t_N]\times\dbR^n$, we define equilibrium strategy function $\mathbbm u^\Pi(\cd\,,\cd)$ on $[t_{N-1},t_N]\times\dbR^n$ as follows (recalling $\psi(\cd)$ from (H2)):
\bel{u(N)}\mathbbm u^\Pi(t,x)=\psi(t_{N-1},t,x,V^\Pi(t,x),V^\Pi_x(t,x),V^\Pi_{xx}(t,x)),\qq(t,x)\in[t_{N-1},t_N]\times\dbR^n.\ee

\it Step 2. \rm Define the {\it representation function} $\Th^{N-1}(\cd\,,\cd)$ on $[t_{N-1},t_N]\times\dbR^n$ by the following equation:
\bel{Th(N-1)}\left\{\2n\ba{ll}
\ns\ds\Th^{N-1}_t(t,x)+\mathbb H\big(t_{N-2},t,x,\mathbbm
u^\Pi(t,x),\Th^{N-1}(t,x),\Th^{N-1}_x(t,x),\Th^{N-1}_{xx}(t,x)\big)=0,\\
\ns\ds\qq\qq\qq\qq\qq\qq\qq\qq\qq\qq\qq\qq\qq\qq(t,x)\in [t_{N-1},t_N]\times\dbR^n,\\
\ns\ds\Th^{N-1}(t_N,x)=h(t_{N-2},x),\qq x\in\dbR^n.\ea\right.\ee
Note that \eqref{Th(N-1)} is different from \eqref{HJB(N)} since $t_{N-1}$ is replaced by $t_{N-2}$.
Consequently, in general, the following fails:
\bel{Th(N-1)=V}\Th^{N-1}(t,x)=V^\Pi(t,x),\qq(t,x)\in[t_{N-1},t_N]\times\dbR^n.\ee
In particular, the following could fail:
\bel{Th(N-1)=V*}\Th^{N-1}(t_{N-1},x)=V^\Pi(t_{N-1},x),\qq x\in\dbR^n.\ee
Now, define the value function $V^\Pi(\cd\,,\cd)$ on $[t_{N-2},t_{N-1})\times\dbR^n$ to be the solution to the following HJB equation:
\bel{HJB(N-1)}\left\{\2n\ba{ll}
\ns\ds V^\Pi_t(t,x)+\inf_{u\in U}\dbH\big(t_{N-2},t,x,u,V^\Pi(t,x),V_x^\Pi(t,x),V^\Pi_{xx}(t,x)\big)=0,\qq
(t,x)\in[t_{N-2},t_{N-1})\times\dbR^n,\\
\ns\ds V^\Pi(t_{N-1}-0,x)=\Th^{N-1}(t_{N-1},x),\qq x\in\dbR^n.\ea\right.\ee
Due to the fact that \eqref{Th(N-1)=V*} fails, $V^\Pi(\cd\,,\cd)$ might have a jump at $t=t_{N-1}$.
Having $V^\Pi(\cd\,,\cd)$ on $[t_{N-2},t_{N-1})\times\dbR^n$, we define equilibrium strategy function $\mathbbm u^\Pi(\cd\,,\cd)$ on $[t_{N-2},t_{N-1})\times\dbR^n$ as follows:
\bel{u(N-1)}\mathbbm u^\Pi(t,x)=\psi(t_{N-2},t,x,V^\Pi(t,x),V^\Pi_x(t,x),V^\Pi_{xx}(t,x)),\qq(t,x)\in[t_{N-2},t_{N-1})
\times\dbR^n.\ee

\it Step 3. \rm Define the {\it representation function} $\Th^{N-2}(\cd\,,\cd)$ on $[t_{N-2},t_N]\times\dbR^n$ by the following equation:
\bel{Th(N-2)}\left\{\2n\ba{ll}
\ns\ds\Th^{N-2}_t(t,x)+\mathbb H\big(t_{N-3},t,x,\mathbbm
u^\Pi(t,x),\Th^{N-2}(t,x),\Th^{N-2}_x(t,x),\Th^{N-2}_{xx}(t,x)\big)=0,\\
\ns\ds\qq\qq\qq\qq\qq\qq\qq\qq\qq\qq\qq\qq\qq\qq(t,x)\in [t_{N-2},t_N]\times\dbR^n,\\
\ns\ds\Th^{N-2}(t_N,x)=h(t_{N-3},x),\qq x\in\dbR^n.\ea\right.\ee
Note that the time interval for \eqref{Th(N-2)} is $[t_{N-2},t_N]$, instead of $[t_{N-2},t_{N-1}]$. Also, unlike \eqref{HJB(N-1)}, $t_{N-3}$ appears, instead of $t_{N-2}$. Thus, in general, the following fails:
\bel{Th(N-2)=V*}\Th^{N-2}(t_{N-2},x)=V^\Pi(t_{N-2},x),\qq x\in\dbR^n.\ee
Now, we define the value function $V^\Pi(\cd\,,\cd)$ on $[t_{N-3},t_{N-2})\times\dbR^n$ to be the solution to the following HJB equation:
\bel{HJB(N-2)}\left\{\2n\ba{ll}
\ns\ds V^\Pi_t(t,x)+\inf_{u\in U}\dbH\big(t_{N-3},t,x,u,V^\Pi(t,x),V_x^\Pi(t,x),V^\Pi_{xx}(t,x)\big)=0,\qq
(t,x)\in[t_{N-3},t_{N-2})\times\dbR^n,\\
\ns\ds V^\Pi(t_{N-2}-0,x)=\Th^{N-2}(t_{N-2},x),\qq x\in\dbR^n.\ea\right.\ee
Since \eqref{Th(N-2)=V*} may fail, $V^\Pi(\cd\,,\cd)$ might have a jump at $t=t_{N-2}$.
Having $V^\Pi(\cd\,,\cd)$ on $[t_{N-3},t_{N-2})\times\dbR^n$, we define equilibrium strategy function $\mathbbm u^\Pi(\cd\,,\cd)$ on $[t_{N-3},t_{N-2})\times\dbR^n$ as follows:
\bel{u(N-2)}\mathbbm u^\Pi(t,x)=\psi(t_{N-3},t,x,V^\Pi(t,x),V^\Pi_x(t,x),V^\Pi_{xx}(t,x)),\qq(t,x)\in[t_{N-3},t_{N-2})
\times\dbR^n.\ee

The rest steps now are clear.

\ms

We write the constructed equilibrium strategy on the whole interval $[0,T]$ for Problem (G$^\Pi$) as follows:
\bel{}\mathbbm u^\Pi(t,x)=\psi\big(\ell^\Pi(t),t,x,V^\Pi(t,x),V^\Pi_x(t,x),V^\Pi_{xx}(t,x)\big),\qq
(t,x)\in[0,T]\times\dbR^n.\ee
Then for any given $x\in\dbR^n$, let $X^\Pi(\cd)$ be the solution to the following SDE:
\bel{}\left\{\2n\ba{ll}
\ns\ds dX^\Pi(s)=b\big(s,X^\Pi(s),\mathbbm u^\Pi(s,X^\Pi(s))\big)ds+\si\big(s,X^\Pi(s),\mathbbm u^\Pi(s,X^\Pi(s))\big)dW(s),\qq s\in [0,T],\\
\ns\ds X^\Pi(0)=x.\ea\right.\ee
Correspondingly, in \eqref{Sec4.3_FBSDE_F_Optimal}, we take
$$\xi_{k-1}=X^\Pi(t_{k-1}),\qq1\les k\les N.$$
Then
\bel{}X^\Pi(s)=\bar X^k\big(s;t_{k-1},\bar X^{k-1}(t_{k-1}),\mathbbm
u^\Pi(\cd)\big),\qq s\in[t_{k-1},t_k),\ k=1,2,\dots,N,\ee
with $\bar X^0(t_0)=x$. For any $t_k\in\Pi\setminus\{t_N\}$, let
$(Y^\Pi(t_k,\cd),Z^\Pi(t_k,\cd))$ be the adapted solution to the
following BSDE:
\bel{}\left\{\2n\ba{ll}
\ns\ds dY^\Pi(t_k,s)=-g\big(t_k,s,X^\Pi(s),\mathbbm
u^\Pi(s,X^\Pi(s)),Y^\Pi(t_k,x),Z^\Pi(t_k,s)\big)ds+Z^\Pi(t_k,s)dW(s),\\
\ns\ds\hskip 9cm s\in [t_k,T],\\
\ns\ds Y^\Pi(t_k,T)=h\big(t_k,X^\Pi(T)\big).\ea\right.\ee
Then from \eqref{Sec4.3-FBSDE(k)}, \eqref{Sec4.3_Y_Theta} and \eqref{Sec4.3_FBSDE_F_Optimal}, we have
\bel{}\left\{\2n\ba{ll}
\ns\ds Y^\Pi(t_{k-1},s)=Y^k(s;t_{k-1},\bar X^{k-1}(t_{k-1}),\mathbbm u^\Pi(\cd)),\\
\ns\ds Z^\Pi(t_{k-1},s)=Z^k(s;t_{k-1},\bar X^{k-1}(t_{k-1}),\mathbbm u^\Pi(\cd)),\ea\right.\qq s\in [t_{k-1},T],\qq k=1,2,\cds,N.\ee

Since, for any $k=1,2,\cds,N$, $\mathbbm u^\Pi(\cd)|_{[t_{k-1},t_k]}$ is an optimal strategy of Problem
(C$_k$), then \eqref{Sec3_Nash} holds true. If we set
$$\bar u^k(\cd)=\mathbbm u^\Pi(\cd\,,X^\Pi(\cd))|_{[t_{k-1},t_k]},\qq 1\les k\les N,$$
then $(\bar u^1(\cd),\bar u^2(\cd),\cds,\bar u^N(\cd))$ can be regarded as a Nash equilibrium of the corresponding $N$-person non-cooperative differential game (see \cite{Yong2012b} for details).

\subsection{The formal limits}

Through an observation on \eqref{Sec4.3_PDE} and \eqref{Sec4.3_HJB}, together with \eqref{Sec4.3_mathbbm_u_F}, we see that $\Th^k(\cd\,,\cd)$ can be extended naturally from $[t_k,t_N]\times\dbR^n$ to $[t_{k-1},t_N]\times\dbR^n$ as follows:
\bel{}\Th^k(t,x)\equiv V^\Pi(t,x),\qq(t,x)\in[t_{k-1},t_k)\times\mathbb R^n,\qq k=1,2,\cds,N.\ee
Then from \eqref{Sec4.3_PDE} and \eqref{Sec4.3_HJB}, we see that extended $\Th^k(\cd\,,\cd)$ satisfies
\bel{Sec4.5_PDE}\left\{\2n\ba{ll}
\ns\ds\Th^k_t(t,x)+\dbH\big(t_{k-1},t,x,\mathbbm u^\Pi(t,x),\Th^k(t,x),\Th^k_x(t,x),\Th^k_{xx}(t,x)
\big)=0,\qq(t,x)\in[t_{k-1},t_N]\times\dbR^n,\\
\ns\ds\Th^k(t_N,x)=h(t_{k-1},x),\qq x\in\dbR^n.\ea\right.\ee
Let us define
\bel{Sec4.5_Notations}\left\{\2n\ba{ll}
\ns\ds\Th^\Pi(\t,t,x)=\sum_{k=1}^N\Th^k(t,x)\mathbbm
1_{[t_{k-1},t_k)}(\t),\qq(\t,t,x)\in D[0,T]\times\dbR^n,\\
\ns\ds h^\Pi(\t,x)=\sum_{k=1}^Nh(t_{k-1},x)\mathbbm
1_{[t_{k-1},t_k)}(\t),\qq(\t,x)\in[0,T]\times\dbR^n,\\
\ns\ds g^\Pi(\t,t,x,u,y,z)=\sum_{k=1}^N
g(t_{k-1},t,x,u,y,z)\mathbbm 1_{[t_{k-1},t_k)}(\t),\\
\ns\ds\hskip 3cm (\t,t,x,u,y,z)\in D[0,T]\times\dbR^n\times
U\times\dbR\times\dbR^{1\times d}.\ea\right.\ee
Then $\Th^\Pi(\cd\,,\cd\,,\cd)$ satisfies the following PDE:
\bel{Sec4.5_PDE_Theta}\left\{\2n\ba{ll}
\ns\ds\Th^\Pi_t(\t,t,x)+\dbH^\Pi\big(\t,t,x,\mathbbm
u^\Pi(t,x),\Th^\Pi(\t,t,x),\Th^\Pi_x(\t,t,x),\Th^\Pi_{xx}(\t,t,x)\big)=0,\\
\ns\ds\hskip 7.5cm (\t,t,x)\in D[0,T]\times\dbR^n,\\
\ns\ds\Th^\Pi(\t,T,x)=h^\Pi(\t,x),\qq(\t,x)\in[0,T]\times\dbR^n,\ea\right.\ee
where
$$\ba{ll}
\ns\ds\dbH^\Pi(\t,t,x,u,\th,p,P)=\tr\big[a(t,x,u)P\big]+\lan b(t,x,u),p\ran
+g^\Pi\big(\t,t,x,u,\th,p^\top\si(t,x,u)\big),\\
\ns\ds\hskip 4.7cm (\t,t,x,u,\th,p,P)\in D[0,T]\times\dbR^n\times U\times\dbR\times\dbR^n\times\dbS^n.\ea$$
Since for any $k=1,2,\cds,N$,
$$\ba{ll}
\ns\ds\mathbbm u^\Pi(t,x)=\psi\big(t_{k-1},t,x,V^\Pi(t,x),V^\Pi_x(t,x),V^\Pi_{xx}(t,x)\big)\\
\ns\ds\qq\qq=\psi\big(t_{k-1},t,x,\Th^k(t,x),\Th^k_x(t,x),\Th^k_{xx}(t,x)\big),\qq(t,x)\in[t_{k-1},t_k)
\times\dbR^n,\ea$$
one has
$$\ba{ll}
\ns\ds\mathbbm u^\Pi(t,x)=\sum_{k=1}^N\psi\big(t_{k-1},t,x,\Th^k(t,x),\Th^k_x(t,x),\Th^k_{xx}(t,x)
\big)\mathbbm1_{[t_{k-1},t_k)}(t)\\
\ns\ds\qq\qq=\psi\(\sum_{k=1}^Nt_{k-1}\mathbbm1_{[t_{k-1},t_k)}(t),t,x,\sum_{k=1}^N\Th^k(t,x)
\mathbbm1_{[t_{k-1},t_k)}(t),\sum_{k=1}^N\Th^k_x(t,x)\mathbbm1_{[t_{k-1},t_k)}(t),\\
\ns\ds\qq\qq\qq\qq\qq\qq\qq\qq\qq\qq\qq\qq\qq\qq\qq\qq\sum_{k=1}^N \Th^k_{xx}(t,x)\mathbbm1_{[t_{k-1},t_k)}(t)\)\\
\ns\ds\qq\qq=\psi\(\ell^\Pi(t),t,x,\Th^\Pi(t,t,x),\Th^\Pi_x(t,t,x),
\Th^\Pi_{xx}(t,t,x)\),\qq(t,x)\in[0,T]\times\dbR^n.\ea$$

Now, we would like to study the behavior when $\|\Pi\| \rightarrow
0$ formally to get the limit equations. In the next section, we will
prove the formal limits can be made rigorous under some conditions.
We introduce the following assumption temporarily.

\ms

{\bf (H3)} There exists a function $\Th(\cd,\cd,\cd)\in C^{0,0,2}(D[0,T]\times\dbR^n)$ such that
$$\ba{ll}
\ns\ds\lim_{\|\Pi\|\to0}\(|\Th^\Pi(\t,t,x)-\Th(\t,t,x)|+|\Th^\Pi_x(\t,t,x)-\Th_x(\t,t,x)|
+|\Th^\Pi_{xx}(\t,t,x)-\Th_{xx}(\t,t,x)|\)=0,\\
\ns\ds\qq\qq\qq\qq\qq\qq\qq\qq\qq\qq\qq\qq\qq(\t,t,x)\in D[0,T]\times\dbR^n,\ea$$
uniformly for $(\t,t,x)$ in any compact sets.

\ms

Under (H3), with the help of (H2), we also have
$$\lim_{\|\Pi\|\to0}|\mathbbm u^\Pi(t,x)-\mathbbm u(t,x)|=0,\qq(t,x)\in[0,T]\times\dbR^n,$$
uniformly for $(t,x)$ in any compact sets, where
\bel{Sec4.5_mathbbm_u}
\mathbbm u(t,x)=\psi\big(t,t,x,\Th(t,t,x),\Th_x(t,t,x),\Th_{xx}(t,t,x)\big),\qq(t,x)\in[0,T]\times\dbR^n.\ee
Therefore, we have
$$\lim_{\|\Pi\|\to0}\| X^\Pi(\cd)-\bar X(\cd)\|_{L^2_\dbF(\O;C([0,T];\dbR^n))}=0,$$
where $\bar X(\cd)$ is the solution to the following SDE:
\bel{}\left\{\2n\ba{ll}
\ns\ds d\bar X(s)=b\big(s,\bar X(s),\mathbbm u(s,\bar X(s))\big)ds+\si\big(s,\bar X(s),
\mathbbm u(s,\bar X(s))\big)dW(s),\qq s\in[0,T],\\
\ns\ds\bar X(0)=x,\ea\right.\ee
and
$$\ba{ll}
\ns\ds L^2_\dbF(\O;C([0,T];\dbR^n))=\Big\{X:[0,T]\times\O\to\dbR^n\bigm|X(\cd)\hb{
is an $\dbF$-progressively measurable process}\\
\ns\ds\qq\qq\qq\qq\qq\qq\qq\hb{with continuous paths such that }\dbE\[\sup_{s\in
[0,T]}|X(s)|^2\]<\infty\Big\}.\ea$$
Moreover, for any $t\in[0,T]$,
$$\lim_{\|\Pi\|\to0}\(\|Y^\Pi(\ell^\Pi(t),\cd)-\bar
Y(t,\cd)\|_{L^2_\dbF(\O;C([0,T];\dbR))}+\|
Z^\Pi(\ell^\Pi(t),\cd)-\bar Z(t,\cd)\|_{L^2_\dbF(0,T;\dbR^{1\times d})}\Big)=0,$$
where $(\bar Y(t,\cd),\bar Z(t,\cd))$ is the solution the following BSDE:
\bel{}\left\{\2n\ba{ll}
\ns\ds d\bar Y(t,s)=-g\big(t,s,\bar X(s),\mathbbm u(s,\bar
X(s)),\bar Y(t,s),\bar Z(t,s)\big)ds+Z(t,s)dW(s),\qq s\in[t,T],\\
\ns\ds\bar Y(t,T)=h\big(t,\bar X(T)\big),\ea\right.\ee
and
$$\ba{ll}
\ns\ds L^2_\dbF(0,T;\dbR^{1\times d})=\Big\{Z:
[0,T]\times\O\to\dbR^{1\times d}\bigm|Z(\cd)
\hb{ is an $\dbF$-progressively measurable process}\\
\ns\ds\qq\qq\qq\qq\qq\qq\hb { such that }\dbE\int_0^T|Z(s)|^2ds<\infty\Big\}.\ea$$
Furthermore, it is clear
$$\ba{ll}
\ns\ds J\big(\ell^\Pi(t),X^\Pi(\ell^\Pi(t)),\mathbbm u^\Pi(\cd)\big)=
Y^\Pi(\ell^\Pi(t),\ell^\Pi(t))=\sum_{k=1}^N Y^k\big(t_{k-1};t_{k-1},\bar
X^{k-1}(t_{k-1}),\mathbbm u^\Pi(\cd)\big)\mathbbm1_{[t_{k-1},t_k)}(t)\\
\ns\ds=\sum_{k=1}^N V^\Pi\big(t_{k-1},\bar X^{k-1}(t_{k-1})\big)\mathbbm 1_{[t_{k-1},t_k)}(t)
=\sum_{k=1}^N\Th^k\big(t_{k-1},\bar X^{k-1}(t_{k-1})\big)\mathbbm1_{[t_{k-1},t_k)}(t)\\
\ns\ds=\Th^\Pi\big(\ell^\Pi(t),\ell^\Pi(t),X^\Pi(\ell^\Pi(t))\big),\qq\qq t\in [0,T].\ea$$
By taking limits, we have
$$J\big(t,\bar X(t),\mathbbm u^\Pi(\cd)\big)=\bar Y(t,t)=\Th\big(t,t,\bar X(t)\big),$$
which coincides with \eqref{Sec3_ValueFucntion}. By Definition \ref{Sec3_Def}, $\mathbbm u(\cd\,,\cd)$ is a time-consistent equilibrium strategy, and
$$V(t,x)\equiv\Th(t,t,x),\qq(t,x)\in[0,T]\times\dbR^n$$
is a time-consistent equilibrium value function of Problem (N).

\ms

Next, we come back to \eqref{Sec4.5_PDE_Theta} and try to derive the
limit equation of \eqref{Sec4.5_PDE_Theta} which is used to
characterize the equilibrium value function $V(\cd\,,\cd)$ or the
more general function $\Th(\cd\,,\cd\,,\cd)$. For this aim, we
rewrite \eqref{Sec4.5_PDE_Theta} in the following integral form:
\bel{Sec4.5_PDE_Theta_Int}\ba{ll}
\ns\ds\Th^\Pi(\t,t,x)=h^\Pi(\t,x)+\int_t^T\dbH^\Pi\big(
\t,s,x,\mathbbm u^\Pi(s,x),\Th^\Pi(\t,s,x),
\Th^\Pi_x(\t,s,x),\Th^\Pi_{xx}(\t,s,x)\big)ds,\\
\ns\ds\qq\qq\qq\qq\qq\qq\qq\qq\qq(\t,t,x)\in D[0,T]\times\dbR^n,\ea\ee
and introduce the following assumption:

\ms

{\bf (H4)} There exists a constant $L>0$ such that
$$|h_\t(\t,x)|+|g_\t(\t,t,x,u,y,z)|\les K, \qq\forall\
(\t,t,x,u,y,z)\in D[0,T]\times\dbR^n\times U\times\dbR\times\dbR^{1\times d}.$$

\ms

Under Assumptions (H1)-(H4), we know
$$\left\{\2n\ba{ll}
\ns\ds\lim_{\|\Pi\|\to0}h^\Pi(\t,x)=h(\t,x),\qq(\t,x)\in[0,T]\times\dbR^n,\\
\ns\ds\lim_{\|\Pi\|\to0}g^\Pi\big(\t,s,x,\mathbbm u^\Pi(s,x),\Theta^\Pi(\t,s,x),\Th^\Pi_x(\t,s,x)^\top
\si(s,x,\mathbbm u^\Pi(s,x))\big)\\
\ns\ds\qq=g\big(\t,s,x,\mathbbm u(s,x),\Th(\t,s,x),\Th_x(\t,s,x)^\top\si(s,x,\mathbbm u(s,x))\big),\qq
(\t,s,x)\in D[0,T]\times\dbR^n,\ea\right.$$
which leads to
$$\ba{ll}
\ns\ds\lim_{\|\Pi\|\to0}\dbH^\Pi\big(\t,s,x,\mathbbm u^\Pi(s,x),\Th^\Pi(\t,s,x), \Th^\Pi_x(\t,s,x),\Th^\Pi_{xx}(\t,s,x)\big)\\
\ns\ds\qq\q=\dbH\big(\t,s,x,\mathbbm u(s,x),\Th(\t,s,x),\Th_x(\t,s,x),\Th_{xx}(\t,s,x)\big),\qq
(\t,s,x)\in D[0,T]\times\dbR^n.\ea$$
Therefore, letting $\|\Pi\|\to0$ in \eqref{Sec4.5_PDE_Theta_Int}, we get the function
$\Th(\cd\,,\cd\,,\cd)$ satisfying the following equation:
\bel{Sec4.5_Equilib_HJB_Int}\ba{ll}
\ns\ds\Th(\t,t,x)=h(\t,x)+\int_t^T\dbH\big(\t,s,x,\mathbbm u(s,x),\Th(\t,s,x),
\Th_x(\t,s,x),\Th_{xx}(\t,s,x)\big)ds,\\
\ns\ds\qq\qq\qq\qq\qq\qq\qq\qq\qq\qq\qq\qq(\t,t,x)\in D[0,T]\times\dbR^n,\ea\ee
or in the differential form:
\bel{Sec4.5_Equilib_HJB}\left\{\2n\ba{ll}
\ns\ds\Th_t(\t,t,x)+\dbH\big(\t,t,x,\mathbbm u(t,x),\Th(\t,t,x),\Th_x(\t,t,x),\Th_{xx}(\t,t,x)\big)=0,\qq (\t,t,x)\in D[0,T]\times\dbR^n,\\
\ns\ds\Th(\t,T,x)=h(\t,x),\qq(\t,x)\in[0,T]\times\dbR^n,\ea\right.\ee
where $\dbH(\cd)$ is defined by \eqref{Sec2_Notations} and $\mathbbm u(\cd\,,\cd)$ is defined by \eqref{Sec4.5_mathbbm_u}. We call \eqref{Sec4.5_Equilib_HJB} or \eqref{Sec4.5_Equilib_HJB_Int} the {\it equilibrium Hamilton-Jacobi-Bellman equation} (equilibrium HJB equation, for short) of Problem (N).

\section{Well-Posedness of the Equilibrium HJB Equation}\label{Sec_Well-posedness_HJB}

In this section, we will present the well-posedness of equation
\eqref{Sec4.5_Equilib_HJB} to some extent. First of all, some
observations on \eqref{Sec4.5_Equilib_HJB} are made as follows:

\ms

$\bullet$ System \eqref{Sec4.5_Equilib_HJB} is a fully nonlinear PDE, but not in a classical form. Note that both $\Th(\t,t,x)$ and $\Th(t,t,x)$ appear in the equation at the same time, where
$\Th(t,t,x)$ is the restriction of $\Th(\t,t,x)$ on $\t=t$, which make the existing theory of fully nonlinear PDEs cannot be applied directly to \eqref{Sec4.5_Equilib_HJB} for its well-posedness.

\ms

$\bullet$ The recursive costs in our time-inconsistent control problem (N) bring some differences from the problem studied in \cite{Yong2012b}, particularly reflecting on the equilibrium HJB equation. If $\Th(t,t,x)$ could be obtained from an independent way, the equilibrium HJB equation developed in \cite{Yong2012b} is in fact a linear PDE with respect to $\Th(\t,t,x)$ where $\t$ could be regarded as a parameter. While in the current situation, it is different now. In fact, \eqref{Sec4.5_Equilib_HJB} will still be a nonlinear one in spite of $\Th(t,t,x)$ is known in advance.

\ms

$\bullet$ Notice the expression: $\mathbbm u(t,x)=\psi\big(t,t,x,\Th(t,t,x),\Th_x(t,t,x),\Th_{xx}(t,t,x)\big)$. From the definition of $\psi$ (see Assumption (H2)), it is clear that the dependence of $\si$ on the control variable $u$ leads to
the appearance of $\Th_{xx}(t,t,x)$ in $\mathbbm u(t,x)$. It turns out that the appearance of $\Th_{xx}(t,t,x)$ will bring some essential difficulties in establishing the well-posedness of equilibrium HJB equation. At the moment, we are not able to overcome the difficulty. We hope to come back in our future publications. In the current paper, having formally derived the general equilibrium HJB equation, we will establish its well-posedness for a special, but still important case. More precisely, we assume that
\bel{Sec5_Assumption_sigma}
\si(t,x,u)=\si(t,x),\qq(t,x,u)\in[0,T]\times\dbR^n\times U,\ee
in the following study.

\ms

Under \eqref{Sec5_Assumption_sigma}, the equilibrium HJB equation \eqref{Sec4.5_Equilib_HJB} reads
\bel{Sec5_Equilib_HJB}\left\{\2n\ba{ll}
\ns\ds\Th_t(\t,t,x)+\tr\big[a(t,x)\Th_{xx}(\t,t,x)\big]+\lan b\big(t,x,\mathbbm u(t,x)),\ \Th_x(\t,t,x)
\ran\\
\ns\ds\qq\qq\qq+g\big(\t,t,x,\mathbbm
u(t,x),\Th(\t,t,x),\Th_x(\t,t,x)^\top
\si(t,x)\big)=0,\qq(\t,t,x)\in D[0,T]\times\dbR^n,\\
\ns\ds\Th(\t,T,x)=h(\t,x),\qq(\t,x)\in[0,T]\times\dbR^n,\ea\right.\ee
where
\bel{u-equ}\mathbbm u(t,x)=\psi\big(t,t,x,\Th(t,t,x),\Th_x(t,t,x)\big),\qq(t,x)\in[0,T]\times\dbR^n.\ee
To avoid heavy notations, we simplify \eqref{Sec5_Equilib_HJB} as follows:
\bel{Sec5_simiplified HJB}\left\{\2n\ba{ll}
\ns\ds\Th_t(\t,t,x)+\tr\big[a(t,x)\Th_{xx}(\t,t,x)\big]+\lan b(t,x,\Th(t,t,x),\Th_x(t,t,x)),\ \Th_x(\t,t,x)\ran\\
\ns\ds\qq\qq+g\big(\t,t,x,\Th(t,t,x),\Th_x(t,t,x),\Th(\t,t,x),\Th_x(\t,t,x)\big)=0,\qq(\t,t,x)\in D[0,T]\times\dbR^n,\\
\ns\ds\Th(\t,T,x)=h(\t,x),\qq(\t,x)\in[0,T]\times\dbR^n,\ea\right.\ee
where
$$\left\{\2n\ba{ll}
\ns\ds b(t,x,\bar\th,\bar p)=b(t,x,\psi(t,t,x,\bar\th,\bar
p)),\qq(t,x,\bar\th,\bar p)\in [0,T]\times\dbR^n\times\dbR\times\dbR^n,\\
\ns\ds g(\t,t,x,\bar\th,\bar p,\th,p)=g(\t,t,x,\psi(t,t,x,\bar\th,\bar
p),\th,p^\top\si(t,x)),\\
\ns\ds\hskip 3.5cm (\t,t,x,\bar\th,\bar p,\th,p)\in
[0,T]\times\dbR^n\times\dbR\times\dbR^n\times\dbR\times\dbR^n.\ea\right.$$

Next, we introduce some spaces. For any $\a\in(0,1)$, we let
\begin{itemize}
\item $C^\a(\dbR^n)$ be the set of all continuous functions:
$\f:\dbR^n\to\dbR$ such that
$$\|\f\|_\a=\|\f\|_0+[\f]_\a<\infty,$$
where
$$\|\f\|_0=\sup_{x\in\dbR^n}|\f(x)|,\qq[\f]_\a=\sup_{x,y\in\dbR^n,x\ne y}{|\f(x)-\f(y)|\over|x-y|^\a},$$

\item $C^{1+\a}(\dbR^n)$ be the set of all continuously differentiable functions: $\f:\dbR^n\to\dbR$
such that
$$\|\f\|_{1+\a}=\|\f\|_0+\|\f_x\|_0+[\f_x]_\a<\infty,$$

\item $C^{2+\a}(\dbR^n)$ be the set of all twice continuously
differential functions: $\f:\dbR^n\to\dbR$ such that
$$\|\f\|_{2+\a}=\|\f\|_0+\|\f_x\|_0+\|\f_{xx}\|_0+[\f_{xx}]_\a<\infty,$$

\item $B([0,T];C^{k+\a}(\dbR^n))$ ($k=0,1,2$) be the set of all measurable
functions $f:[0,T]\times\dbR^n\to\dbR$ such that for any $t\in[0,T]$, $f(t,\cd)\in C^{k+\a}(\dbR^n)$ and
$$\|f(\cd\,,\cd)\|_{B([0,T];C^{k+\a}(\dbR^n))}=\sup_{t\in
[0,T]}\|f(t,\cd)\|_{k+\a}<\infty,$$

\item $C([0,T];C^{k+\a}(\dbR^n))$ ($k=0,1,2$) be the subset of $B([0,T];C^{k+\a}(\dbR^n))$ consisting of all continuous functions.
\end{itemize}

\ms

We also need the following assumption.

\ms

{\bf (A)} The maps $a:[0,T]\times\dbR^n\to\dbS^n,\ b:[0,T]\times\dbR^n\times\dbR\times
\dbR^n\to\dbR^n,\ g:D[0,T]\times\dbR^n\times\dbR\times\dbR^n\times\dbR\times\dbR^n\to\dbR$, and
$h:[0,T]\times\dbR^n\to\dbR$ are continuous and bounded. Moreover, there exists a constant $L>0$ such that
$$\ba{ll}
\ns\ds|a_x(t,x)|+|b_x(t,x,\bar\th,\bar p)|+|g_x(\t,t,x,\bar\th,\bar p,\th,p)|+|b_{\bar \th}(t,x,\bar\th,\bar p)|+|g_{\bar\th}(\t,t,x,\bar\th,\bar p,\th,p)|+|b_{\bar p}(t,x,\bar\th,\bar p)|\\
\ns\ds+|g_{\bar p}(\t,t,x,\bar\th,\bar p,\th,p)|+|g_\th(\t,t,x,\bar\th,\bar p,\th,p)|+|g_p(\t,t,x,\bar\th,\bar p,\th,p)|+|h_x(\t,x)|+|h_{\t x}(\t,x)|\les L,\\
\ns\ds\hskip 7cm (\t,t,x,\bar\th,\bar p,\th,p)\in D[0,T]\times
\dbR^n\times\dbR\times\dbR^n\times\dbR\times\dbR^n.\ea$$
Furthermore, $a(t,x)^{-1}$ exists for all $(t,x)\in [0,T]\times
\dbR^n$ and there exist constants $\l_0,\ \l_1>0$ such that
$$\l_0 I\les a(t,x)^{-1}\les\l_1I,\qq(t,x)\in[0,T]\times\dbR^n.$$

For any $v(\cd\,,\cd)\in C([0,T];C^1(\dbR^n))$, we consider
the following semi-linear PDE parameterized by $\t\in [0,T]$:
\bel{Sec5_semi-linear PDE}\left\{\2n\ba{ll}
\ns\ds\Th_t(\t,t,x)+\tr\big[a(t,x)\Th_{xx}(\t,t,x)\big]+\lan b(t,x,v(t,x),v_x(t,x)),\ \Th_x(\t,t,x)\ran\\
\ns\ds\hskip 1cm+g\big(\t,t,x,v(t,x),v_x(t,x),\Th(\t,t,x),\Th_x(\t,t,x)\big)=0,\qq(t,x)\in[0,T]\times \dbR^n,\\
\ns\ds\Th(\t,T,x)=h(\t,x),\qq x\in\dbR^n,\ea\right.\ee
which admits a unique solution by the classical theory. In order to
derive the well-posedness of \eqref{Sec5_simiplified HJB}, we will
utilize the above semi-linear PDE to establish a contraction mapping
from $v(\cd\,,\cd)$ to $V(\cd\,,\cd)$, where $V(t,x)
=\Th(t,t,x)$, $(t,x)\in [0,T]\times\dbR^n$. The details will
be presented in the following theorem.

\begin{theorem} \sl Under Assumption (A), equation \eqref{Sec5_simiplified HJB} admits a
unique solution.
\end{theorem}

\it Proof. \rm By the fundamental solution theory of parabolic PDE (see Friedman
\cite{Friedman1964}), for any $v(\cd\,,\cd)\in C([0,T];C^1(\dbR^n))$, the solution of \eqref{Sec5_semi-linear PDE} can be expressed as
\bel{Sec5_solution}\ba{ll}
\ns\ds\Th(\t,t,x)=\int_{\dbR^n}\G(t,x;T,y)h(\t,y)dy+\int_t^T\int_{\dbR^n}\G(t,x;s,y)\[\lan b(s,y, v(s,y),v_x(s,y)),\Th_x(\t,s,y)\ran\\
\ns\ds\qq\qq+g(\t,s,y,v(s,y),v_x(s,y),\Th(\t,s,y),\Th_x(\t,s,y))\]dyds,\qq(t,x)\in[\t,T]\times\dbR^n,\ea\ee
where $\G(t,x;s,y)$ called the fundamental solution is defined
on $[0,T]\times\dbR^n\times[0,T]\times\dbR^n$ with $t<s$ and
$$\G(t,x;s,y)={1\over(4\pi(s-t))^{n\over2}(\det[a(s,y)])^{1\over2}}e^{-\lan a(s,y)^{-1}(x-y),x-y\ran \over4(s-t)}.$$
Some direct calculation leads to the following estimates:
\bel{Sec5_Estimates}\left\{\2n\ba{ll}
\ns\ds|\G(t,x;s,y)|\les K{e^{-{\l|x-y|^2\over4(s-t)}}\over(s-t)^{n\over2}} ,\\
\ns\ds|\G_x(t,x;s,y)|\les K{e^{-{\l|x-y|^2\over4(s-t)}}\over(s-t)^{n+1\over2}},\ea\right.\qq \l<\l_0,\ee
where $K>0$ is a constant which can be different from line to line. Moreover,
\bel{Sec5_Gamma_y}\G_y(t,x;s,y)=-\G_x(t,x;s,y)-\G(t,x;s,y)\rho(t,x;s,y),\ee
where
$$\left\{\2n\ba{ll}
\ns\ds\rho(t,x;s,y)={(\det[a(s,y)])_y\over2\det[a(s,y)]}+{\lan[a(s,y)^{-1}]_y(x-y),\ x-y\ran\over4(s-t)},\\
\ns\ds\lan[a(s,y)^{-1}]_y(x-y),\ x-y\ran=\left(\begin{array}{ccc}
 \lan[a(s,y)^{-1}]_{y_1}(x-y),\ x-y\ran \\
  \lan [a(s,y)^{-1}]_{y_2}(x-y),\ x-y\ran \\
  \vdots\\
 \lan[a(s,y)^{-1}]_{y_n}(x-y),\ x-y\ran\end{array}  \right).\ea\right.$$
It is easy to check, under (A),
\bel{}|\rho(t,x;s,y)|\les K\(1+{|x-y|^2\over s-t}\).\ee

We now split the rest of the proof into several steps.

\ms

\it Step 1. \rm We prove  $\Th_x(\t,t,x)$ is bounded, i.e.,
$$|\Th_x(\t,t,x)|\les K(1+\|h(\t,\cd)\|_{C^1(\dbR^n)}),\qq(t,x)\in[\t,T]\times\dbR^n.$$
Combined \eqref{Sec5_Gamma_y} with the method of integration by parts, we get
\bel{}\ba{ll}
\ns\ds\int_{\dbR^n}\G_x(t,x;T,y)h(\t,y)dy=-\int_{\dbR^n}\G_y(t,x;T,y)h(\t,y)dy-\int_{\dbR^n}\G(t,x; T,y)\rho(t,x;T,y)h(\t,y)dy\\
\ns\ds=\int_{\dbR^n}\G(t,x;T,y)h_y(\t,y)dy-\int_{\dbR^n}\G(t,x;T,y)\rho(t,x;T,y)h(\t,y)dy,\ea\ee
then
\bel{Sec5_x-y-2}\ba{ll}
\ns\ds\Th_x(\t,t,x)=\int_{\dbR^n}\G(t,x;T,y)h_y(\t,y)dy-\int_{\dbR^n}\G(t,x;T,y)\rho(t,x;T,y)h(\t,y)dy\\
\ns\ds\qq\qq+\int_t^T\int_{\dbR^n}\G_x(t,x;s,y)\[\lan b(s,y,v(s,y),v_x(s,y)),\Th_x(\t,s,y)\ran\\
\ns\ds\hskip3.5cm+g(\t,s,y,v(s,y),v_x(s,y),\Th(\t,s,y),\Th_x(\t,s,y))\]dyds.\ea\ee
Therefore, from Estimates \eqref{Sec5_Estimates} and Assumption (A), we get
\bel{}\ba{ll}
\ns\ds|\Th_x(\t,t,x)|\les\int_{\dbR^n}K{e^{-{\l|x-y|^2\over4(T-t)}}\over(T-t)^{n\over2}}\[|h_y(\t ,y)|+\(1+{|x-y|^2\over T-t}\)|h(\t,y)|\]dy\\
\ns\ds+\int_t^T\int_{\dbR^n}K{e^{-{\l|x-y|^2\over4(s-t)}}\over(s-t)^{n+1\over2}}\big(|\Th_x(\t,s, y)|+1\big)dyds\\
\ns\ds\les K\(1+\|h(\t,\cd)\|_{C^1(\dbR^n)}\)+\int_t^T\int_{\dbR^n}K {e^{-{\l|x-y|^2\over4(s-t)}}\over(s-t)^{n+1\over2}}|\Th_x(\t,s,y)|dyds.\ea\ee
By Gronwall's inequality, we obtain
\bel{Sec5_boundedness}
|\Th_x(\t,t,x)|\les K\(1+\|h(\t,\cd)\|_{C^1(\dbR^n)}\),\qq(t,x)\in[\t,T]\times\dbR^n.\ee

\it Step 2. \rm  For any $v^1(\cd\,,\cd),v^2(\cd\,,\cd)\in C([0,T];C^1(\dbR^n))$, let $\Th^1(\cd\,,\cd\,,\cd)$ and $\Th^2(\cd\,,\cd\,,\cd)$ be the corresponding solutions to \eqref{Sec5_solution}. We want to prove
\bel{}\ba{ll}
\ns\ds\Th^1(\t,t,x)-\Th^2(\t,t,x)|+|\Th^1_x(\t,t,x)-\Th^2_x(\t,t,x)|\les K
(T-t)^{1\over2}\|v^1(\cd\,,\cd)-v^2(\cd\,,\cd)\|_{C([\t,T];C^1(\dbR^n))}\ea\ee
holds true for any $(t,x)\in [\t,T]\times\dbR^n$. To this end, for any $(\t,t,x)\in D[0,T]\times\dbR^n$, we denote
$$\left\{\2n\ba{ll}
\ns\ds\D v(t,x)=v^1(t,x)-v^2(t,x),\\
\ns\ds\D\Th(\t,t,x)=\Th^1(\t,t,x)-\Th^2(\t,t,x),\\
\ns\ds\D b(t,x)=b(t,x,v^1(t,x),v_x^1(t,x))-b(t,x,v^2(t,x),v_x^2(t,x)),\\
\ns\ds\D g(\t,t,x)=g(\t,t,x,v^1(t,x),v_x^1(t,x),\Th^1(\t,t,x),\Th^1_x(\t,t,x))\\
\ns\ds\hskip2cm-g(\t,t,x,v^2(t,x),v_x^2(t,x),\Th^2(\t,t,x),\Th^2_x(\t,t,x)).\ea\right.$$
Then, we have
\bel{}\ba{ll}
\ns\ds|\D\Th(\t,t,x)|=\Big|\int_t^T\int_{\dbR^n}\G(t,x;s,y)\[\lan\D b(s,y),\Th^1_x(\t,s,y)\ran\\
\ns\ds\qq\qq\qq\qq+\lan b(s,y,v^2(s,y),v_x^2(s,y)),\D\Th_x(\t,s,y)\ran+\D g(\t,s,y)\]dyds\Big|.\ea\ee
By Estimates \eqref{Sec5_Estimates}, \eqref{Sec5_boundedness}, and
Assumption (A),
\bel{Sec5_Theta^1-Theta^2}\ba{ll}
\ns\ds|\D\Th(\t,t,x)|\les\int_t^T\int_{\dbR^n}K{e^{-{\l|x-y|^2\over4(s-t)}}\over(s-t)^{n\over2}}
\[\(|\D v(s,y)|+|\D v_x(s,y)|\)\(1+|\Th^1_x(\t,s,y)|\)\\
\ns\ds\qq+|\D\Th(\t,s,y)|+|\D\Th_x(\t,s,y)|\]dyds\\
\ns\ds\les K(T-t)\(1+\|h(\t,\cd)\|_{C^1(\dbR^n)}\)\|\D v(\cd\,,\cd)\|_{C([\t,T];C^1(\dbR^n))}\\
\ns\ds\qq+\int_t^T\int_{\dbR^n}K{e^{-{\l|x-y|^2\over4(s-t)}}\over(s-t)^{n\over2}}\[
|\D\Th(\t,s,y)|+|\D\Th_x(\t,s,y)|\]dyds.\ea\ee
Similarly,
\bel{Sec5_Theta^1x-Theta^2x}\ba{ll}
\ns\ds|\D\Th_x(\t,t,x)|\les K(T-t)^{1\over2}\(1+\|h(\t,\cd)\|_{C^1(\dbR^n)}\)\|\D v(\cd\,,\cd)\|_{C([\t,T];C^1(\dbR^n))}\\
\ns\ds\qq+\int_t^T\int_{\dbR^n}K{e^{-{\l|x-y|^2\over4(s-t)}}\over(s-t)^{n+1\over2}}\[
|\D\Th(\t,s,y)|+|\D\Th_x(\t,s,y)|\]dyds.\ea\ee
Then, by combining \eqref{Sec5_Theta^1-Theta^2} and
\eqref{Sec5_Theta^1x-Theta^2x}, we have
\bel{Sec5_sum of Theta^1-Theta^2}\ba{ll}
\ns\ds|\D\Th(\t,t,x)|+|\D\Th_x(\t,t,x)|\les K(T-t)^{1\over2}\(1+\|h(\t,\cd)\|_{C^1(\dbR^n)}\)
\|\D v(\cd\,,\cd)\|_{C([\t,T];C^1(\dbR^n))}\\
\ns\ds\qq\qq\qq\qq\qq+\int_t^T\int_{\dbR^n}K{e^{-{\l|x-y|^2\over4(s-t)}}\over(s-t)^{n+1\over2}}\[
|\D\Th(\t,s,y)|+|\D\Th_x(\t,s,y)|\]dyds.\ea\ee
Gronwall's inequality works again to yield
\bel{Sec5_norm of Theta^1-Theta^2 [tau,T]}\|\D\Th(\t,\cd\,,\cd)\|_{C([\t,T];C^1(\dbR^n))}
\les K(T-t)^{1\over2}\(1+\|h(\t,\cd)\|_{C^1(\dbR^n)}\)\|\D v(\cd\,,\cd)\|_{C([\t,T];C^1(\dbR^n))}.\ee
We note that the constant $K>0$ appearing in the above is independent of $(\t,t)\in D[0,T]$.

\ms

\it Step 3. \rm Denoting $V^i(t,x)=\Th^i(t,t,x)$
($i=1,2$) and taking $\t=t$ in \eqref{Sec5_norm of
Theta^1-Theta^2 [tau,T]}, we get
\bel{Sec5_norm of V^1-V^2 [t,T]}\ba{ll}
\ns\ds\|V^1(\cd\,,\cd)-V^2(\cd\,,\cd)\|_{C([t,T];C^1(\dbR^n))}\\
\ns\ds\les K(T-t)^{1\over2}\(1+\|h(\cd\,,\cd)\|_{B([0,T];C^1(\dbR^n))}\)\|v^1(\cd\,,\cd)-v^2(\cd\,,\cd)\|_{C([t,T];
C^1(\dbR^n))}.\ea\ee
Obviously, a contraction mapping $v(\cd\,,\cd)\mapsto V(\cd\,,\cd)$ on $C([T-\d,T];C^1(\dbR^n))$ is obtained by choosing $\d>0$ small enough. Accordingly, this map has a unique fixed point on $[T-\d,T]$. Furthermore, similar estimates on $[T-2\d,T-\d], [T-3\d,T-2\d],\cds$, till to $[0,\d]$ are derived as above so that the map $v(\cd\,,\cd)\mapsto V(\cd\,,\cd)$ admits a unique fixed point
on the whole space $C([0,T];C^1(\dbR^n))$. Therefore, we get the well-posedness of the following integro-differential equation:
\bel{Sec5_sloved HJB}\ba{ll}
\ns\ds\Th(\t,t,x)\1n=\2n\int_{\dbR^n}\2n\G(t,x;T,y)h(\t,y)dy\1n+\1n\int_t^T\3n\int_{\dbR^n}\2n\G(t,x;s,y)
\[\lan b(s,y,\Th(s,s,y),\Th_x(s,s,y)),\Th_x(\t,s,y)\ran\\
\ns\ds\qq\qq\qq+g(\t,s,y,\Th(s,s,y),\Th_x(s,s,y),\Th(\t,s,y),\Th_x(\t,s,y))\]dyds,\q(\t,t,x)\1n\in\1n D[0,T]\1n\times\1n\dbR^n.\ea\ee
Finally, by the classical theory of PDE, and the regularity of the above expression, we get $\Th(\t,t,x)$ is $C^{2+\a}$ in $x$, $C^{1+{\a\over2}}$ in $t$ for some $\a\in(0,1)$, and PDE \eqref{Sec5_simiplified HJB} is satisfied. \endpf

\ms

In the previous section, we introduced Assumption (H3) to provide the convergence of $\Th^\Pi(\cd\,,\cd\,,\cd)$ in the space $C^{0,0,2}(D[0,T]\times\mathbb R^n)$, which guarantees the
existences of an equilibrium strategy and a corresponding
equilibrium value function of Problem (N). However, as $\sigma$ does
not depend on $u$ (see Condition \eqref{Sec5_Assumption_sigma}),
Assumption (H3) is reduced to the convergence in the space
$C^{0,0,1}(D[0,T]\times\dbR^n)$, i.e., (H3) is
replaced by the following:

\ms

{\bf (H3$'$).} There exists a function $\Th(\cd\,,\cd\,,\cd)\in C^{0,0,1}(D[0,T]\times\dbR^n)$
such that
\bel{Sec5_convergence assumption}
\lim_{\|\Pi\|\to0}\(|\Th^\Pi(\t,t,x)-\Th(\t,t,x)|+|\Th^\Pi_x(\t,t,x)-\Th_x(\t,t,x)|\)=0,\qq(\t,t,x)\in D[0,T]\times\dbR^n,\ee
uniformly for $(\t,t,x)$ in any compact sets.

\ms

To make the study rigorously, we shall prove the expected convergence \eqref{Sec5_convergence assumption} holds true exactly for the family $\Th^\Pi(\cd\,,\cd\,,\cd)$ constructed in the previous section.
We still require all the involved functions are bound and continuously differentiable up to a needed order with bounded derivatives. Consequently, the uniform Lipschitz continuous of $\t\mapsto(h(\t,x),h_x(\t,x),g(\t,t,x,y,z,u))$ follows from (H4) and (A). Therefore
$$\ba{ll}
\ns\ds|h^\Pi(\t,x)-h(\t,x)|+|h^\Pi_x(\t,x)-h_x(\t,x)|+|g^\Pi(\t,t,x,u,y,z)-g(\t,t,x,u,y,z)|\les
K \|\Pi\|,\\
\ns\ds\qq\qq\qq\qq\qq\qq\qq(\t,t,x,u,y,z)\in D[0,T]\times\dbR^n\times U\times
\dbR\times\dbR^{1\times d},\ea$$
where $g^\Pi$ and $h^\Pi$ are defined by \eqref{Sec4.5_Notations}.

\ms

\it Proof of \eqref{Sec5_convergence assumption}. \rm Under Condition \eqref{Sec5_Assumption_sigma},
\eqref{Sec4.5_PDE_Theta} reads
\bel{Sec5 Theta^Pi}\left\{\2n\ba{ll}
\ns\ds\Th^\Pi_t(\t,t,x)+\tr\[a(t,x)\Th^\Pi_{xx}(\t,t,x)\]+\lan b\big(t,x,\psi(\ell^\Pi(t),t,x,\Th^\Pi(t,t,x),\Th^\Pi_x(t,t,x))\big),
\ \Th^\Pi_x(\t,t,x)\ran\\
\ns\ds\qq\qq+g^\Pi\big(\t,t,x,\psi(\ell^\Pi(t),t,x,\Th^\Pi(t,t,x),\Th^\Pi_x(t,t,x)),\Th^\Pi(\t,t,x),
\Th^\Pi_x(\t,t,x)^\top\si(t,x)\big)=0,\\
\ns\ds\hskip10cm (\t,t,x)\in D[0,T]\times\dbR^n,\\
\ns\ds\Th^\Pi(\t,T,x)=h^\Pi(\t,x),\qq(\t,x)\in[0,T]\times\dbR^n.\ea\right.\ee
By \eqref{Sec5_solution}, we have
\bel{Sec5_Theta^Pi-Theta}\ba{ll}
\ns\ds\D\Th^\Pi(\t,t,x)=\int_{\dbR^n}\G(t,x;T,y)\D h^\Pi(\t,y)dy+\int_t^T\int_{\dbR^n}\G(t,x;s,y)
\[\lan\D b^\Pi(s,y),\ \Th^\Pi_x(\t,s,y)\ran\\
\ns\ds\qq\qq+\lan b(s,y,\psi(s,s,y,\Th(s,s,y),\Th_x(s,s,y))), \
\D\Th^\Pi_x(\t,s,y)\ran+\D g^\Pi(s,y)\]dyds,\ea\ee
where we denote
$$\left\{\2n\ba{ll}
\ns\ds\D\Th^\Pi(\t,s,y)=\Th^\Pi(\t,s,y)-\Th(\t,s,y),\\
\ns\ds\D h^\Pi(\t,y)=h^\Pi(\t,y)-h(\t,y),\\
\ns\ds\D b^\Pi(s,y)=b\big(s,y,\psi(\ell^\Pi(s),s,y,\Th^\Pi(s,s,y),\Th_x^\Pi(s,s,y))\big)- b\big(s,y,\psi(s,s,y,\Th(s,s,y),\Th_x(s,s,y))\big),\\
\ns\ds\D g^\Pi(s,y)=g^\Pi(\t,s,y,\psi(\ell^\Pi(s),s,y,\Th^\Pi(s,s,y),\Th_x^\Pi(s,s,y)),\Th^\Pi(\t,s,y),
\Th^\Pi_x(\t,s,y)^\top\si(s,y))\\
\ns\ds\hskip2cm-g(\t,s,y,\psi(s,s,y,\Th(s,s,y),\Th_x(s,s,y)),\Th(\t,s,y),\Th_x(\t,s,y)^\top\si(s,y)).\ea
\right.$$
Similarly,
\bel{Sec5_Theta^Pi_x-Theta^Pi}\ba{ll}
\ns\ds\D\Th^\Pi_x(\t,t,x)=\int_{\dbR^n}\G(t,x;T,y)\D h^\Pi_y(\t,y)dy-\int_{\dbR^n}\G(t,x;T,y)\rho(t,x,T,y) \D h^\Pi(\t,y)dy\\
\ns\ds\qq\qq\qq+\int_t^T\int_{\dbR^n}\G_x(t,x;s,y)\[\lan\D b^\Pi(s,y),\ \Th^\Pi_x(\t,s,y)\ran\\
\ns\ds\qq\qq\qq+\lan b\big(s,y,\psi(s,s,y,\Th(s,s,y),\Th_x(s,s,y))\big),\ \D\Th^\Pi_x(\t,s,y)\ran+\D g^\Pi(s,y)\]dyds.\ea\ee
Therefore, from \eqref{Sec5_Theta^Pi-Theta} and \eqref{Sec5_Theta^Pi_x-Theta^Pi}, one obtains
$$\ba{ll}
\ns\ds|\D\Th^\Pi(\t,t,x)|+|\D\Th^\Pi_x(\t,t,x)|\\
\ns\ds\les\int_{\dbR^n}{Ke^{-{\l|x-y|^2\over4(T-t)}}\over(T-t)^{n\over2}}\cd
\|\Pi\|dy+\int_{\dbR^n}{Ke^{-{\l|x-y|^2\over4(T-t)}}\over(T-t)^{n\over2}}\(1+{|x-y|^2\over4(T-t)}\)\cd
\|\Pi\|dy\\
\ns\ds\qq+\int_t^T\int_{\dbR^n}{Ke^{-{\l|x-y|^2\over4(s-t)}}\over(s-t)^{n+1\over2}}\[\|\Pi\|
+|\D\Th^\Pi(s,s,y)|+|\D\Th^\Pi_x(s,s,y)|+|\D\Th^\Pi(\t,s,y)|+|\D\Th^\Pi_x(\t,s,y)|\]dyds\\
\ns\ds\les K\|\Pi\|+\int_t^T\int_{\dbR^n}{Ke^{-{\l|x-y|^2\over4(s-t)}}\over(s-t)^{n+1\over2}}
\[|\D\Th^\Pi(s,s,y)|+|\D\Th^\Pi_x(s,s,y)|+|\D\Th^\Pi(\t,s,y)|+|\D\Th^\Pi_x(\t,s,y)|\]dyds.\ea$$
Then,
$$\ba{ll}
\ns\ds\sup_{\t\in[0,t]}\(|\D\Th^\Pi(\t,t,x)|+|\D\Th^\Pi_x(\t,t,x)|\)\les K\|\Pi\|\\
\ns\ds\qq\qq+\int_t^T\int_{\dbR^n}{Ke^{-{\l|x-y|^2\over4(s-t)}}\over(s-t)^{n+1\over2}}
\sup_{\t\in[0,s]}\(|\D\Th^\Pi(\t,s,y)|+|\D\Th^\Pi_x(\t,s,y)|\)dyds.\ea$$
Applying Gronwall's inequality leads to
$$\sup_{\t\in[0,t]}\(|\D\Th^\Pi(\t,t,x)|+|\D\Th^\Pi_x(\t,t,x)|\)\les K\|\Pi\|,\qq t\in[0,T].$$
We have proved the expected convergence \eqref{Sec5_convergence assumption}.\endpf

\section*{Acknowledgement}

This work was carried out during the stay of Qingmeng Wei and
Zhiyong Yu at University of Central Florida, USA. They would like to
thank the hospitality of Department of Mathematics, and the
financial support from China Scholarship Council.

\end{document}